\documentclass[12pt]{article}
\usepackage{amssymb,amsmath,amsthm,latexsym}

\title{A projection construction for semifields%
\thanks{This research was supported in part by
 NSA grant H98230-10-1-0159  and 
 NSF grant  DMS~0753640.}}

\author{J\"urgen Bierbrauer%
\\
Department of Mathematical Sciences\\
Michigan Technological University\\
Houghton, MI 49931 \\ \\
William M. Kantor%
\\
Department of Mathematics\\
University of Oregon \\
Eugene, OR 97403  }
\begin{document}
\maketitle

\newtheorem{Theorem}[equation]{Theorem}
\newtheorem{theorem}[equation]{Theorem}
\newtheorem{Proposition}[equation]{Proposition}
\newtheorem{proposition}[equation]{Proposition}
\newtheorem{lemma}[equation]{Lemma}
\newtheorem{Lemma}[equation]{Lemma}
\newtheorem{definition}[equation]{Definition}
\newtheorem{Definition}[equation]{Definition}
\newtheorem{Corollary}[equation]{Corollary}
\newtheorem{corollary}[equation]{Corollary}
\newtheorem{conjecture}[equation]{Conjecture}
\newtheorem{fact}[equation]{\hspace{-3pt}\rm}
\newtheorem{Remark}[equation]{Remark}
\newtheorem{remark}[equation]{Remark} 

\newtheorem{Remarks}[equation]{Remarks}
\newtheorem{hypothesis}[equation]{Hypothesis}
\newtheorem{coro}[equation]{Corollary}
\newtheorem{prop}[equation]{Proposition}
\newtheorem{exam}[equation]{Example}
\newtheorem{Example}[equation]{Example}
\newtheorem{Examples}[equation]{Examples}
\newtheorem{notation}[equation]{Notation}
\newtheorem{convention}[equation]{Convention}

\def\para#1{\medskip\noindent{\bf #1}~}

\numberwithin{equation}{section}

\def\nz{\mathbb{N}}
\def\gz{\mathbb{Z}}
\def\mz{\mathbb{M}}
\def\wz{\mathbb{W}}
\def\rz{\mathbb{R}}
\def\ef{\mathbb{F}}
\def\CC{\mathbb{C}}
\def\o{\omega}
\def\p{\overline{\omega}}
\def\e{\epsilon}
\def\a{\alpha}
\def\b{\beta}
\def\g{\gamma}
\def\d{\delta}
\def\l{\lambda}
\def\s{\sigma}
\def\t{\tau}
\def\O{\Omega}
\def\sm{\setminus}
\def\bsl{\backslash}
\def\la{\longrightarrow}
\def\arr{\rightarrow}
\def\ov{\overline}
\def\Aut{\rm Aut}
\newcommand{\D}{\displaystyle}
\newcommand{\T}{\textstyle}

\def\mar#1{{\marginpar{\small#1}}}

\def\mod{{\rm mod }}
\abstract{Simple constructions are given for finite semifields that include  
as special cases both old semifields and recently constructed semifields.}

\section{Introduction}
\label{Introduction}

A finite {\em semifield} is a finite non-associative algebra 
$(F,+,*)$ (i.e.,  not necessarily associative algebra) with identity element 
$1$ such that $x*y=0\Rightarrow x=0$ or $y=0$.  
Thus, $(F,+)$ is an abelian group, and both distributive laws hold. 
A finite {\em presemifield} satisfies the same conditions except for the 
existence of an identity element.  Semifields first arose in the study of algebras 
resembling fields \cite{Di1}, and soon afterwards were used to construct 
finite nondesarguesian projective planes \cite{VW}.
A fundamental result of Albert \cite{Albert isomorphism theorem} states that 
two semifields (or presemifields) $(F,+,*)$ and $(F',+,*')$ determine 
isomorphic planes if and only if the algebras are {\em isotopic}: 
there are additive group isomorphisms $A, B, C\colon F\to F'$ such that 
$A(x)*'B(y)=C(x*y)$ for all $x,y\in F$.  Each presemifield is isotopic 
to a semifield (cf. Section~\ref{Isotopes}). There are many surveys
of finite semifields (e.g., \cite[Ch. V]{De} and \cite{Pingree}).

If $F=\oplus_iU_i$ is a direct decomposition of the group $(F,+)$, 
with corresponding  projection maps $E_i$, then $x*y=\sum_iE_i(x*y)$.  
We will slightly generalize this trivial observation and reverse this process 
(Proposition~\ref{noncommprojprodprop}) in order to obtain both new views 
of well-known semifields and new semifields.  
In Section~\ref{basicsection} we collect some general information concerning 
some presemifields, the corresponding semifields
and their nuclei. In Section~\ref{twofieldssection} it is shown that 
the Dickson semifields \cite{Di1}, the  Hughes-Kleinfeld semifields \cite{HK} 
and some of the Knuth semifields \cite{KnuthJAlg} can be constructed via
Proposition~\ref{noncommprojprodprop} using as ingredients semifields  
isotopic to fields.

Our new constructions are based on   
(isotopes of) fields and the {\em twisted fields} of Albert \cite{Alberttwist}.   
Our presemifields and some of their properties  are obtained in 
Sections~\ref{constrAsection}-\ref {oddBCnucleisection}, 
and  are summarized in the following

\begin{Theorem} 
\label{gentheorem}
There are presemifields of order $p^{2m}$  with $p$ a prime  
of the following sorts $($where $\s =p^s$
with $1\leq s<m ;$   $v\in L=\ef_{p^{m}};$ 
$  l,n, N, R\in L^*;$ $  \b\in F^*=\ef_{p^{2m}}^*)$$:$
\begin{itemize}
\item[\rm(a)] $A(p,m,s,l,\mu ) $ 
with  $-l \notin L^{\s-1} $ and
  $  \ov{\mu}/\mu $   not in the subgroup of $\,F^*{} ^{p^m-1} $ of index  
$ \gcd(p^m+1,\s+1).$ 
\item[\rm(b)] $B(p,m,s, l,n, N)$ 
 with $p$ odd$,$ where    $-l\notin L^{\s -1} $ and 
 the polynomial
 $t^{\s +1}+(1-1/N)t+(1/n-1)/N$  has no root in $L$.

\item[\rm(c)] $C(p,m,s, l,R)$ 
 with $p$ odd$,$ where    $-l\notin L^{\s -1} $ and 
  $R\notin L^{\s +1}.$  
  \end{itemize}
For odd $p$  the nuclei are described in 
{\rm Theorems~\ref{BClrnuketheorem},  \ref{BspecialNmtheorem},  
 \ref{not N squared},}  and   
{\rm  Propositions~\ref{special case of of nuclei}, \ref{last case}.}
\end{Theorem}

Some of these semifields are
  isotopic  to commutative ones, but  most are not when $p$ is odd.
 Section~\ref{commsection} handles this isotopism question for odd $p$, using a criterion due to 
  Ganley~\cite[Theorem~4]{Ganley72}. We use a semifield  associated with a given
 presemifield (see Proposition~\ref{pretosemigenprop}) to extend the Ganley criterion to
 presemifields (Corollary~\ref{preGanleyprop}).
A simple, direct proof of the Ganley criterion  is given in  Theorem~\ref{Ganleyprop}.

The right  nucleus $\nz_r $,  middle  nucleus $\nz_m$ and  left nucleus 
$\nz_l$  of a semifield $F$  are the sets
 of all $r, m$ or $l\in F$ behaving as follows  for all $x,y\in F$. 
 $${ 
 \begin{array}{c|c|cc}
 r\in \nz_r  & m\in \nz_m  &   l\in \nz_l \\
 \hline
 \!(x*y)*r \! = \! x*(y*r)\! & \!(x*m)*y \! = \! x*(m*y) 
\!& \!(l*x)*y \! = \! l*(x*y)\!
 \end{array}
 }
 $$
\noindent These nuclei are significant for the planes 
(it is an old and elementary result about projective planes   
\cite[Theorem~8.2]{HP} that   the nuclei  correspond to three natural 
homology groups -- important types of   collineation  groups). 
In particular, isotopic 
semifields have isomorphic nuclei.
We emphasize that nuclei are only dealt with in odd characteristic.  
In that case, in all parts of Theorem~\ref{gentheorem} the left and right nuclei are identical,
of order $p^{\gcd(s,m)}.$ Whenever we are able to determine the middle nucleus, it turns
out to be at most 
a quadratic extension of the left (right) nucleus.
Moreover, the method used in Section~\ref{Generalization for all odd $q$} 
in odd characteristic 
in order to obtain the $B$- and $C$-families from the $A$-family undoubtedly 
works 
in characteristic $2,$ but the calculations become tedious and so have 
been avoided. The efforts  used in Section~\ref{oddBCnucleisection}
to study nuclei would again be more detailed in characteristic 2.
 
Several  instances of our constructions have already appeared.  
The families of commutative semifields constructed in
Budaghyan-Helleseth \cite{BuHe11} are special cases of
$A(p,m,s,1/N,\mu )$ in odd characteristic.
Their first family is  isotopic to members of the $C$-family
(see Theorem~\ref{Ccommtheorem} and Proposition~\ref{BuHeCprop});
the members of their second family   not in the
first are isotopic to $B(p,m,s,1/\o^{\s-1},n,n^{(\s-1)/2})$,
where $n\in L$ is a non-square 
(see Proposition~\ref{BuHeBprop}).
\par
Planar functions possess at least two different analogues in characteristic $2$:
commutative semifields and APN functions.
The latter form a core part of the theory of cryptographic $S$-boxes.
The projection method works also in this context, and 
the APN hexanomials of Budaghyan-Carlet~\cite{BC08} are  characteristic 
$2$ analogues of a special case of the $A$-family.

\section{A projection construction}
\label{introsection}
 All algebras will be finite.
 
\begin{Definition}
\label{noncommprojproddef} 
Let  $(F,+,*_i),$   $i=1,\dots ,n,$
  be  non-associative algebras with the same additive group $F,$
  and   let $A_i\subset F$ be additive subgroups. Then the set 
$\lbrace (*_i,A_i)\rbrace$ is called {\em compatible} if
$x*_iy\in A_i$ for $i=1,\dots ,n$ imply that $x$ or $y=0.$
\end{Definition}
The following is evident:
 
 \begin{Proposition}
\label{noncommprojprodprop}
Let $F= \oplus_iU_i .$
Suppose that $(F,+,*_i)$ are  non-associative algebras$,$ $i=1,\dots ,n,$ 
such that the set
$\lbrace (*_i,A_i)\rbrace$ is compatible and 
$\dim _{\ef_{p}}A_i +\dim _{\ef_{p}} U_i=\dim _{\ef_{p}} F .$
If $f_i\colon F\to U_i$ is any $\ef_p$-linear map such that  $\ker(f_i)=A_i,$ 
then
\begin{equation}
\label{general fi}
x\circ y=\sum_if_i(x*_iy)
\end{equation}
defines a presemifield on $F.$ 

\end{Proposition}
 
Conversely:
\begin{Proposition}
\label{getallprop}

Any presemifield $(F,+,*)$ can be constructed in many ways using 
{\rm Proposition~\ref{noncommprojprodprop}}.
 \end{Proposition}
\begin{proof}
Use an arbitrary direct sum decomposition $F=\oplus_{i=1}^nU_i$,
$U_i\ne0$, 
with $f_i$ the projections  to the  $U_i$.
Define $x*_iy= f_i( x*y)\in U_i.$ 
If $A_i=\oplus_{j\not=i}U_j$  then the set $\lbrace (*_i,A_i)\rbrace$ 
is compatible, and 
Proposition~\ref{noncommprojprodprop}  yields $(F,+,*)$.
\end{proof}

The preceding proof shows that the construction in 
Proposition~\ref{noncommprojprodprop}
yields all presemifields, and each one in ridiculously many ways. 
We will make this construction more
interesting by using ingredients $x*_iy$ that are 
sufficiently  elementary.

Although the proof of Proposition~\ref{noncommprojprodprop} is trivial, 
it has several interesting
instances. It was first used for odd characteristic  commutative semifields 
in \cite{JBprojsemi}, 
where the
generalization to the non-commutative case was already alluded to.
In the applications we will concentrate on the case when $n=2,r=2m$ and 
$r_1=r_2=m.$
\emph{The advantage of the notion of compatibility is that the maps $f_i$
 are essentially ignored.}  In the next section we will discuss the isotopy 
of some of the presemifields (\ref{general fi}) for different choices of 
these maps.

\section{Basic considerations}
\label{basicsection}
\subsection{Projections}
Assume that 
$F=\ef_{q^{2m}} \supset L=\ef_{q^{ m}} ,$ $ x\to\ov{x}$ is the involutory 
automorphism 
of $F$ and $T\colon F\to L$ is the trace map. 

Let $\o\in F^*$ with $T(\o)=0,$ where
 $\o=1$ in  characteristic $2$. 

Assume that $(*_1,u_1L)$ and $(*_2,u_2L) $ are compatible, 
where $u_i\in F$ 
(Definition~\ref{noncommprojproddef}). 
We will consider the corresponding presemifield $(F,*)$ 
of Proposition~\ref{noncommprojprodprop}
up to isotopy, using the   $L$-linear mappings $f_i(x)=T((\o /u_i)x)z_i$
on $F$, $i=1,2$,
where $z_i\in F^*$ and $z_2/z_1\notin L.$

Up to isotopy  we can choose
$z_1=1$ and $z_2=z$ an arbitrarily fixed element not in $L.$
The corresponding presemifield product   therefore has  the form
$x*y=T\big((\o /u_1)(x*_1y) \big)+T \big((\o /u_2)(x*_2y) \big)z.$ 
Let $*'_i$ be defined by $x*'_iy=(x*_iy)/u_i.$ Then 
$(*'_1,L)$ is compatible with $(*'_2,L)$ and we have
$x*y=T \big(\o (x*'_1y) \big)+T \big(\o (x*'_2y) \big)z.$
Replacing the products $x*_iy$ by multiples $\b_i (x*_iy)$ we can 
therefore assume
that the compatibility condition has   $A_1=A_2=L.$  
Using any $z\notin L$, (\ref{general fi}) has the form
\begin{equation}
\label{two term definition}
x*y=T \big(\o (x*_1y) \big)+T \big(\o (x*_2y) \big)z.
\end{equation}
Of course, different choices for $z$ produce isotopic presemifields.
For example, in odd characteristic we can choose $z=1/\o ,$ 
which corresponds to $f_2(x)=x-\ov{x}$ and
$$x*y=\frac{1}{2}\big [\o (x*_1y- \ov{x*_1y})+(x*_2y-\ov{x*_2y}) \big]. $$

\para{Neighbors.} 
In the setting of  Definition~\ref{noncommprojproddef}  we call 
$n$-tuples $\{(*_i,A_i)\}$  and 
$\{(*'_i,A_i)\}$ 
{\em neighbors}  
if  $x*_iy-x*'_iy\in A_i$ for all $x,y\in F.$ 
Then $\{(*_i,A_i)\}$ compatible if and only if 
$\{(*'_i,A_i)\}$ is.
\emph{The resulting presemifield products} (\ref{general fi})
\emph{are then identical.}
 Thus, apparently different ingredients $\{(*_i,A_i)\}$ can produce the same 
semifields.

For example, consider the field product $xy$ in odd characteristic and a 
Dickson-Knuth product
$x\underset{\!\! \scriptscriptstyle DK\!\!}\circ y=ac+b^{\s}d^{\tau}\mu +
(ad+bc)\omega$ where $x=a+b\omega ,y=c+d\omega$
and $a,b,c,d\in L.$ Then $x\underset{\!\! 
\scriptscriptstyle DK\!\!}\circ y-xy\in L$ for all $x,y\in L.$ 
Whenever $(xy,L)$ is used as one member of a compatible pair it can be 
replaced by $(x\underset{\!\! \scriptscriptstyle DK\!\!}\circ y,L)$
without changing the resulting presemifield.

\subsection{Isotopes}
\label{Isotopes} 
 {\em Any presemifield
$(F,+,*)$  is isotopic to a  semifield}.  The traditional approach
\cite[p.~957]{Kaplansky}   is 
to fix \emph{any}
$0\ne e\in F$
  and define $\star$ by 
  \begin{equation}
  \label{isotope identity}
  (x*e) \star(e*y)=x*y
    \end{equation}
     for all $
x,y\in F$.  Then $ (F,+, \star )$ is a semifield with
identity element
$e*e$, and 
is obviously isotopic to $(F,+,*)$.
We will use a different isotope with 
 the property  that $e$, not $e*e$, is the identity element.   
Define $\b ,\g\colon F\to F$ by 
\begin{equation}
\label{beta gamma}
\mbox{$e* \b (x)=x   $ 
\ and \  $\g (x)*e=e*x   $}
\end{equation}
(so that $\g(e)=e$, $\g(\b(x ))*e=x$ and $\b(e*y)=y$), 
and 
  \begin{equation}
  \label{isotope identity e}
  x\circ y := \b\big (\g (x)*y \big).
    \end{equation}
Then {\em $ (F,+,\circ)$ is a semifield 
 isotopic to $(F,+,*)$ 
with identity element $e$}, since
$
x\circ e = \b \big(\g (x)*e\big)= \b (e*x)=x$ and 
$
e\circ y = \b \big(\g (e)*y \big)= \b (e*y)=y$
for any $x,y\in F$. If  $(F,+,* )$ is commutative then so are   
$ (F,+, \star)$
and    $ (F,+, \circ)$ (since $\g(x)=x)$).


The following observations build  on $ (F,+,\circ)$ to obtain information 
concerning the nuclei of any semifield isotopic to a given presemifield.
(The {\em center} of a semifield is the set of elements in  all three nuclei 
that commute with all elements.)

\begin{Proposition}
\label{pretosemigenprop}  
Consider fields $F_1\subseteq F_ 2\subseteq F$  and a presemifield 
$(F,+,*)$. 
Assume  that
\begin{itemize}
\item [\rm (i)] $(ax)*y=x*(ay)=a(x*y)$ for all $x,y\in F$,  $a\in F_1,$  and
\item  [\rm(ii)] $(xk)*y=x*(ky)$ for all $x,y\in F, k\in F_2.$
\end{itemize}
If $\b ,\g$ and $x\circ y$ are defined as above using $e=1,$ then   
$(F,+,\circ )$ is a~semifield with identity element $1,$ center containing 
$F_1$  and   middle
nucleus containing  $F_2,$ such that
$k\circ y=yk$ for all $y\in F$, $k\in F_2.$
\
 
 \end{Proposition}

\begin{proof} 
 By (i,ii),  $\b $ and $\g$ are $F_1 $-linear,
 so that 
  $x\circ a =ax= a\circ x$ for all $x\in F,a\in F_1.$ 
  Then  (i) and (\ref{isotope identity e}) imply that  $F_1$ is a subfield of 
 the center  of   $F$. 

If  $k\in F_2,$  $ y\in F,$   we obtain the 
final assertion of the proposition:
$$k\circ y= \b \big (\g (k)*y \big) = \b ([1k]*y)= \b (1*[ky] )=ky $$
 by (ii)  and the definitions of $\b, \g$ and $\circ$.
 Finally,   we claim that  $k $ is in the middle nucleus
of $(F,+,\circ )$. 
For, if also  $x\in F$, then,
again  by   (ii)  and definition,  
$$x\circ (k\circ y)=x\circ (yk)=\b [\g (x)*(yk) ]=\b [\big(k\g (x) \big)*y],$$
 while 
$(x\circ k)\circ y=\b [\g (x\circ k)*y].$ \
Finally, for the same reasons,  $\g(x \circ k) = 
\g \big( \b(\g(x )  *[k1] ) \big)=\g \big(\b([k\g(x)]*1) \big) =k\g(x),$  
as required.
\end{proof}

\begin{Proposition}
\label{pretosemigenpropNUCLEI}

Let  $F$ be a field and $(F,+,*)$ a presemifield$,$  and define 
$\b ,\g$ and $x\circ y$ as above  using $e=1$.

\begin{itemize}
\item  [\rm(a)] Any $x\in F$ satisfying $(\g (x)y)*z=\g (x)(y*z)$ for all 
$y,z\in F $   is in the left nucleus
of $(F,\circ , +).$ 
\item  [\rm(b)] 
Any $z\in F$ satisfying  $x*(yz)=z(x*y)$ for all $x,y \in F$   is in the right
 nucleus of $(F,\circ , +)$ and satisfies $y\circ z=yz$ for all $y\in F.$
\end{itemize}
\end{Proposition}

\begin{proof}
(a) We have  $\g (x\circ y)=\g (x)\g (y)$ for all $y $ 
since
$\g (x\circ y)*1=1*(x\circ y)=\g(x)*y$ by the definitions 
of $\circ $ and $\g$, while the definition  
of $\g$ and two applications of our hypothesis concerning $x$ yield
$[\g (x)\g (y)]*1= \g(x)[\g(y)*1]=\g (x)(1*y)=[\g(x)1]*y$.

By the definition of $\circ$ and hypothesis, if $z\in F$ then 
  $1*[(x\circ y)\circ z]=\g(x\circ y)*z = [\g(x)\g(y)] *z =\g(x)[\g(y)*z]
  =\g(x)[\g(y)*z]=
 [\g(x)1]*(y\circ z) =1*[x\circ (y \circ z)]$, 
so that $x$ is in the left nucleus.

\smallskip
(b) By   definition  and hypothesis,  if $x,y\in F$ then 
$ 1*(y\circ z ) =\g (y)*(1z)= z[\g (y)*1]= z(1*y)= 
1*(yz )  $ 
and hence
$1*[(x\circ y)\circ z]
= \g(x\circ y)*[1z] 
=z[1*(x\circ y)]
=z[\g(x)*y]
=\g(x)*(yz)
=1*[x\circ (yz)]
=1*[x\circ (y\circ z)] $, 
so that $z$ is in the right nucleus.
\end{proof}
\begin{Corollary}
\label{pretosemicor} 
Let $(F,+,*_i)$  and the compatible  set 
$\lbrace (*_i,A_i)\rbrace$ be as in {\rm Definition~\ref{noncommprojproddef}}  
for   $i=1,\dots ,n$. 
Use the  presemifield multiplication    
{\rm(\ref{general fi}).}
Let $ F_1\subseteq F_2 $ be subfields of $F$ 
behaving as in  
{\rm Proposition~\ref{pretosemigenprop}(i),(ii),} 
 for each    $(F,+,*_i),$  
where  the $f_i$ are $F_1$-linear.
Then any semifield  isotopic to $(F,+,*)$
has center of order at least $|F_1|$ and middle nucleus of order at least 
$|F_2|$.
\end{Corollary}
\begin{proof} As the $*_i$ satisfy Proposition~\ref{pretosemigenprop}(i),(ii),  the $F_1$-linearity
of the $f_i$ implies that the presemifield multiplication
{\rm(\ref{general fi})} also satisfies Proposition~\ref{pretosemigenprop}(i),(ii). The preceding proposition
yields our claim.
\end{proof}

\subsection{Isotopes of  commutative semifields}
The next result is due to 
 Ganley~\cite[Theorem~4]{Ganley72}, who used  a very different type of proof.

\begin{Theorem}
\label{Ganleyprop}
A semifield $(F,+,\circ )$ is isotopic to a commutative presemifield if and only if
there is some $w\not=0$ such that $(w\circ x)\circ y=(w\circ y)\circ x$ for all $x,y\in F.$
\end{Theorem}
\begin{proof} If  $w$ exists, then 
 $x\cdot y=(w\circ x)\circ y$ defines a commutative presemifield isotopic   to $(F,+,\circ )$.
\par

Conversely, if $(F,+,\circ )$ is a semifield isotopic to a commutative presemifield $(F,+,\cdot )$, then there are (additive) permutations $\a, \b, \g$ of $F$ such that $\g(x\cdot y)=\a(x)\circ \b(y)$,
and hence $\a(x)\circ \b(y)=\a(y)\circ \b(x)$, for all $x, y\in F$.  Let $\b(u)=1$ (the identity element of $(F,+,\circ )$) 
and   $w=\a(u)$.
Then $\a(x) =\a(x)\circ \b(u)=\a(u)\circ \b(x)=w\circ \b(x)$, and hence 
$(w\circ \b(x))\circ \b(y)=
(w\circ \b(y))\circ \b(x)$  for all $x,y\in F.$
\end{proof}

Of course, all commutative presemifields are isotopic to commutative ones
by (\ref{isotope identity})  or  (\ref{isotope identity e}).
The preceding proof did not use the additive group at all:  it really provides  a simple criterion for a loop to be isotopic to a commutative loop (and a similar statement holds for (\ref{isotope identity})  
and   (\ref{isotope identity e}).)

\begin{Corollary}
\label{preGanleyprop}
Let $(F,+,*)$ be a presemifield. Let the maps $\b,\g$ in 
{\rm (\ref{beta gamma})} define  the associated semifield {\rm (\ref{isotope identity e}).} 
The function $\a(x)=\g(\b(x))$ is defined by $\a(x)*1=x.$
Then $(F,+,*)$ is isotopic to a commutative semifield if and
only if there is some $v\not=0$ such that $\a (v*x)*y=\a (v*y)*x$ for all $x,y.$
\end{Corollary}
\begin{proof}
If  $v\not=0$ exists, then $x\cdot y=\a (v*x)*y$ 
defines an isotopic  commutative presemifield. 

Conversely, if  $(F,+,*)$ is isotopic to a commutative presemifield then so is the semifield $(F,+,\circ )$  in (\ref{isotope identity e}). By Theorem~\ref{Ganleyprop},
there is some $w\not=0$ such that $(w\circ x)\circ y=(w\circ y)\circ x$ for all $x,y.$
Let $v=\g (w).$ Then, by  (\ref{isotope identity e}), $(w\circ x)\circ y=\b (\a(v*x)*y)$  also equals  $(w\circ y)\circ x=\b (\a(v*y)*x)$.
\end{proof}

We will use Corollary~\ref{preGanleyprop} to show that many of our new
semifields are isotopic to commutative ones and many are not.

\section{Field multiplications as ingredients}
\label{twofieldssection}

It is natural to ask which presemifields can be constructed using as 
ingredients only
products   isotopic to field multiplication.
In this section we briefly discuss examples showing that it is possible 
to apply Proposition~\ref{noncommprojprodprop} to
isotopic copies of the same semifield product 
and obtain non-isotopic semifields.

\para{The Hughes-Kleinfeld semifields.}
Consider  $F=\ef_{q^4}$ in odd characteristic, $T\colon F\to L= \ef_{q^2}$ 
the trace, and $\o\in F^*$ of trace $0 $
so that  $n= \o^2$ is a non-square in $L.$
  Write elements of $F$ in the form $x=a+b\o  $ with $a,b\in L.$
If  $l\in L\sm \ef_q, $ define the Hughes-Kleinfeld semifield 
$(F,+, \underset{\scriptscriptstyle HK}{\circ})$   \cite{HK}  by
$$(a+b\o )\underset{\!\! \scriptscriptstyle HK\!\!}{\circ}(c+d\o )=
ac+lb^qd +(a^qd+bc)\o .$$
If also $y=c+d\o$, let $x*_1y=ac+lb^qd$ and $x*_2y=(a^qd+bc)\o$ 
be the projections of $\underset{\scriptscriptstyle HK}{\circ}$ to $L$ and
$L\o ,$ respectively. 
Define   $\ef_q$-linear 
invertible mappings  $\a_i\colon F\to F$ by
$\a_1(a+b\o )=a+(l/n )b^q\o $
and
$\a_2(a+b\o)=a^q+b\o .$
Let $x*'_iy=\a_i(x)y.$ 
Then $*'_1$ is an $L\o -$neighbor of $*_1$, and $*'_2$ is an $L-$neighbor of
$*_2.$ It follows that $\underset{\scriptscriptstyle HK}{\circ}$ 
can be described by Proposition~\ref{noncommprojprodprop}
using as ingredients $x*'_1y$ and $x*'_2y$, both of which are isotopic 
to the field multiplication, since 
$$x\underset{\!\! \scriptscriptstyle HK\!\!}{\circ}y=
(1/2)T\big (\a_1(x)y) \big)+(1/2n )T \big(\o\a_2(x)y \big)\o .$$

\para{The Dickson semifields.}
Let $F=\ef_{q^{2m}}\supset L=\ef_{q^{m}} $ in odd characteristic and 
$\o\in F$ with 
$\o^2=n$ a non-square in $L$.  Let $\a (a+b\o )=a+b^{\s}\o ,$ where $\s$ is a 
power of $q.$
Then $(\o\a (x)\a (y),L)$ is compatible with $(xy,L)$, and the resulting 
product
$$x*y=(1/2)T(\a (x)\a (y))+(1/2)T(xy/\o )\o$$
corresponds to a Dickson semifield  \cite{Di1}.
 
\para{Some Knuth semifields.}
Next we use  Proposition~\ref{noncommprojprodprop}  
to construct one of Knuth's families of semifields $(F,*)$
using as ingredients multiplications  
isotopic to field multiplication.
Let $L=\ef_q\subset F=\ef_{q^2}$; 
choose $f,g\in L$  and  a power  $\s$  of $q $  such that 
  \begin{equation}
  \label{Knuth condition}
\mbox{  $t^{\s+1}+tg-f\not=0$ for all $t\in L.$}
   \end{equation} 
If $1,\l$ is a basis of $F$ over $ L $
with $N(\l )=-f$,  then 
\begin{equation}
\label{Knuth}
(a+b\l)*(c+d\l)=(ac+b^{1/\s}df)+(bc+a^{\s}d+bdg)\l 
\end{equation}
(with $a,b,c,d\in L$) 
produces a semifield $(F,+,*)$ \cite[Sec.~7.3]{KnuthJAlg}.

Define three invertible $L$-linear maps $  F\to F $   by 
$\a_1(a+b\l )=a+b^{1/\s}\l ,$ $ \a_2(a+b\l )=a^{\s}+\l  $  and
$\b_2(a+b\l )=a+b(g-T(\l ))+b\l .$

\begin{Proposition}
\label{Knuth4linkprop}

The pairs $(\ov{\l}\a_1(x)y,L)$ and  $(\a_2(x)\b_2(y),L) $ are   
compatible$,$  and a semifield
obtained via  {\rm Proposition~\ref{noncommprojprodprop}}
is isotopic to {\rm(\ref{Knuth}).}
\end{Proposition}
\proof
This is a routine calculation. We have $\l^2=T(\l )\l +f.$ 
If $x=a+g\l $  and  $y=c+d\l ,$
then

\hspace{20pt} 
\vspace{-2pt}
$
\begin{array}{rlllll}
\a_1(x)y &\hspace{-6pt} \in\hspace{-6pt}&
\vspace{3pt} ac+b^{1/\s}df+L\l
\\a_2(x)\b_2(y)&\hspace{-6pt} \in\hspace{-6pt}&  
(a^{\s}d+bc+bd(g-T(\l ))+bdT(\l ))\l +L.  
\hfill \qed
\end{array}
$
\medskip

Thus,  it is possible for the ingredients to be isotopic to commutative 
presemifields
while the new presemifield is not isotopic to a commutative one.

\section{The semifields 
 $A(p,m,s,l, \mu  )$}
\label{constrAsection}
 In this section we will define a  family of examples obtained 
using 
Definition~\ref{noncommprojproddef} and discuss its properties.
In the next section we will describe variations and partial generalizations. 
 
Let $L=\ef_{p^m}\subset F=\ef_{p^{2m}}$,  $l\in L^*$,
$\s =p^s$  with $1< s <m$.
Let  the operation $x\underset{\scriptscriptstyle A}{\circ}y=
x^{}y^{\s}+lx^{\s}y^{} $ on $F$ 
imitate a twisted field (here the subscript $A$ stands for Albert).
Let $\mu \in F^* $.
We apply Definition~\ref{noncommprojproddef} with 
$$
\mbox{$A_1=A_2=L$,
\  $x*_1y= xy,$   ${x*_2y}= 
\mu (\ov{ x}\underset{\scriptscriptstyle A}{\circ}y)$.}
$$
Thus, by 
(\ref{two term definition})  our binary operation is  
\begin{equation}
\label{Aproduct}
x*y=T(\o { x} y  )+T \big(\o \mu(\ov{ x}
\underset{\scriptscriptstyle A}{\circ}y) \big)\o .
\end{equation}
\begin{Theorem} [The semifields $A(p,m,s,l, \mu  )$]
\label{Atheorem}
The pairs $( xy,L)$ and 
 $(\mu  (\ov{x}\underset{\scriptscriptstyle A}{\circ}y),L) $ are compatible
 if and only if  
 \begin{itemize} 
\item [\rm(i)]$-l \notin L^{\s-1} $ and
\item [\rm(ii)] $ \mu ^{p^m-1}=\ov{\mu}/\mu $  is not contained in the 
subgroup of $\,\ef^*{} ^{p^m+1} $ of index  $ \gcd(p^m+1,\s+1)$.
\end{itemize}

\end{Theorem}

\begin{proof}Let $x, y\ne0$.~Then 
 $T(\o xy)=0 \! \iff \! xy\in L^* \! \iff\!  {y=a\ov{x}\ne0}$ for some 
$a\in L^*,$ in which case 
$$
  x*_2y= \mu \big(\ov{x} \underset{\scriptscriptstyle A}{\circ}
  (a\ov{x}) \big)
= \mu\ov{x}^{\s +1} (a^{\s } +l a).
$$
The last factor is in $L,$  and
  (i) is equivalent to
this factor  being nonzero for $a\not=0.$ 
Moreover, 
 (ii) is equivalent to 
 $\mu\notin L^*F ^*{}^{\s +1}$
 and hence to $\mu\ov{x}^{\s +1}\notin L$
for  all  $x\in F^*.$ 
\end{proof}

We emphasize the elementary nature of this proof using a  
convenient factorization.  
Note that  (i) is precisely the condition needed for 
$(L,\underset{\scriptscriptstyle A}{\circ})$ to be a twisted field, so that 
  $(F,*)$ can be viewed as a ``quadratic extension of a twisted field'' 
(compare Lemma~\ref{Albertsublemma}).
   The gcd in (ii) needs to be $>1$.  Hence, if $\s=1$ then $q$ must be odd; 
but in this case our presemifield is isotopic to a field.

   There is a more general version of (\ref{Aproduct}) using  
   $x\underset{\scriptscriptstyle A}{\circ}y=x^{\s'}y^{\s}+lx^{\s}y^{\s'} $
   with $\s'\in \Aut F$, and 
$ x*y=T(\o \mu_1{ x} y  )+T \big(\o \mu_2(\ov{ x}
\underset{\scriptscriptstyle A}{\circ}y) \big)z  $ for  $\mu_i\in F^* $ 
such that $ \mu_1^{\s -\s'}\in L.$
We leave it to the reader to verify that this apparently more general product 
yields a presemifield isotopic to   one of those in  the theorem.

\begin{Corollary}
\label{firstfamnukecor}
There is a semifield isotopic to $A(p,m,s, l, \mu  )$ such that $\ef_q$ is 
in the center while
the field $\lbrace z\mid z\in F,~ \ov{z}=z^{\s}=z^{1/\s}\rbrace $
is in the middle nucleus. 
\end{Corollary}
\begin{proof} This follows from Corollary~\ref{pretosemicor}.
\end{proof}

\begin{Lemma}
\label{famAisotoplemma}
Up to isotopy of 
$A(p,m,s, l, \mu  ),$ $l$ can be replaced by an arbitrary   member of the coset
$lL^*{}^{\s-1}$ and $\mu $ can be replaced by an arbitrary member of $\mu L^*{}.$
\end{Lemma}
\begin{proof}
Use $x'=\a x,$ $ y'=y/\a  $ as  new variables to map 
$$l\to l/N(\a)^{\s-1}, \ \mu\to \mu\a^{\s}/\ov{\a}^{ } .$$

The last statement corresponds to using $z'=kz,$ $ k\in L^*$.
\end{proof}

\begin{Lemma}
\label{Albertsublemma} 
There is a presemifield   $(F,\star)$ isotopic to   $A(p,m,s,l, \mu)$ such that
the subfield $L$ of $F$ is closed under $\star $ and  
$(L, +,\star)$ is isotopic to 
$(L,+,\underset{\scriptscriptstyle A}{\circ})$.
\end{Lemma}
\begin{proof} 
 An isotopic presemifield product is
$x \star y=T(\o \mu (x\underset{\scriptscriptstyle A}{\circ} y))+
T(\o\ov{x}y)\o $. 
If $x,y\in L $ then
$x \star y=T(\o \mu )({x}\underset{\scriptscriptstyle A}{\circ}y)$, 
so that 
$(L,+,\star)$ and $(L,+,\underset{\scriptscriptstyle A}{\circ}\,)$ 
are isotopic.  
\end{proof}

It follows from the preceding lemma that $A(p,m,s,l, \mu )$ can be isotopic
to a field only if the corresponding twisted field is. 
This is the case if and only if $\s^2$ is the identity on $L.$ 
  (In the latter case the twisted field 
$(L,+,\underset{\scriptscriptstyle A}{\circ})$ is isotopic to 
$(L,+,{\circ '})$, $x\circ'y:=xy+l(xy)^\s$, and hence to a field.)
 
 In Section~\ref{oddBCnucleisection}  we will see that
Corollary~\ref{firstfamnukecor} and Lemma~\ref{Albertsublemma}
provide information on the nuclei
as the nuclei of  twisted fields are known \cite[Theorem~1]{Alberttwist}.
The conditions in Theorem~\ref{Atheorem} can be met if and only if
both $\gcd(q^m-1,\s-1)$ and $\gcd(q^m+1,\s+1)$ are not $1.$ 
These requirements always hold in odd characteristic, so the construction is 
more widely applicable in
that case.  If $p=2$ then the conditions in Theorem~\ref{Atheorem} 
can be met
if $\log_2q^m$ is an odd integer divisible by $\log_2s.$ 

\para{Characteristic $2$.}
We briefly consider   $A(2,m,s,l, \mu )$ in Theorem~\ref{Atheorem}. If
$Q=2^m $ then  $F=\ef_{Q^2}\supset L=\ef_{Q}.$  Let $m=2u$ for odd $u $  
and choose $\sigma =4.$ Then $\gcd(Q-1,\s -1)=3$ and 
our requirement on $l$ is  that 
$l\in L\sm L^3.$ It follows from Lemma~\ref{famAisotoplemma} 
that the choice of $l$ 
does not affect the isotopy class.
As $u$ is odd we have $\gcd(Q+1,\s +1)=5.$ 
It follows from Lemma~\ref{famAisotoplemma} that
we may assume that
$\mu $ has norm $1$ and is not a fifth power in the norm $1$ group.
Observe that $\ef_{16}\subset F$ and $25$ does
not divide $Q^2-1$ unless $u=5$ \mod $10.$ If $u$ is not $5$ \mod $10$ 
we can choose
$\mu $ as an element of order $5.$
\par
The smallest case, $A(2,2,2,\o, \mu )$ of order $16,$ is a field.
Here $\o$ is an element of order $3,$ the field $\ef_{16}$ has been 
constructed as
$\ef_2(\a )$ where $\a^4=\a+1$ and $\mu =\a^3,\o =\a^5.$ 
The presemifield multiplication can be chosen as 
$x*y=T(\mu (\ov{xy}+\o xy))+T(xy) \mu $ where $T\colon \ef_{16}\to\ef_4$ 
is the
trace. Proposition~\ref{pretosemigenprop} shows that this is a field.
\par
The smallest interesting case is $A(2,6,2,l, \mu )$, of order $2^{12}$,
where $\mu =\a^3\in\ef_{16}$ as
before and $l\in L=\ef_{64}$ has order $9$ and $l^3=\o .$ The presemifield
product is
$$(x,y)=T(\mu (\ov{x}y^4+l\ov{x}^4y))+T(xy) \mu .$$
Proposition~\ref{pretosemigenprop} applies. The corresponding semifield has
left and right nucleus $\ef_4$ and middle nucleus $\ef_{16}.$
We do not know whether  $A(2,6,2,l, \mu )$ is isotopic to a twisted field.

\section{The semifields $B(p,m,s,l,n,N)$   and  
$C(p,m,s,l,R)$}
\label{Generalization for all odd $q$}

For odd $p$ we now rewrite the product  
(\ref{Aproduct}) entirely in terms of $L$. This will allow us 
to generalize the presemifields 
$A(p,m,s,l, \mu )$.

 Let $F$ and $L$ be as usual with  
  $p$ odd.
Let $T(\o)=0\ne \o $, so that $n:=\o^2$ is a nonsquare in $L^*$,
 and $u\in F $ satisfies $T(\o u)=0 $ if and only if  $  u\in L$.
Identify $F$ with $L^2$ via   $(a, b)=a+b\o$ and 
\begin{equation}\label{field}
(a,b) (c,d)=(a c+nb d,ad+bc)
\end{equation}
 whenever $a,b,c,d\in L$.  
Then $\ov{(a,b)} =(a,-b)$. 
For $x=a+b\o =(a,b)$ we also write $a=\pi_\mathfrak{R}(x), b=\pi_\mathfrak{I}(x)\,$
(we think of these as the ``real'' and ``imaginary'' parts of $x$).
If $\s=p^s$ then 
$(a, b)^\s=(a+b\o)^\s  = 
 (a^\s ,N b^\s )$ where $N :=\o^{\s-1}\in L$.
 

  \medskip
  
  {\noindent \bf Motivation.}
Let $x=(a,b), y=(c,d).$
Up to a constant factor, the presemifield multiplication 
in  Theorem~\ref{Atheorem} is
$$x\star y=(\pi_\mathfrak{I}\big (\mu (\ov{x}\underset{\scriptscriptstyle A}{\circ}y),\pi_\mathfrak{I}(xy) \big).$$
By Theorem~\ref{Atheorem}(ii),  $\bar\mu\ne \mu$ and hence
$\mu\notin L$.
 Then we may assume that  $\mu =(nv,1)$ 
by the second part of  Lemma~\ref{famAisotoplemma}.

A straightforward calculation shows that  
the product is
\begin{equation}
\label{genformequ}
(a,b) \star(c,d) =(h(a,b,c,d),ad+bc),
\end{equation}
where
 \begin{eqnarray}
\label{A'''' mult}
\label{B mult}
\begin{array}{llllll}
h(a,b,c,d)  &\hspace{-6pt}:= &\hspace{-6pt} \{[a^{{}}c^{\s }-n N  b^{{}}  d^{\s }] 
     &\hspace{-6pt}+ &\hspace{-6pt}l[a^{\s}c^{{} }-  nN  b^{\s}  d^{{} }]\}
     \vspace{4pt}
     \\
&&  \hspace{-27pt} +\:nv \{[Na^{{}}  d^{\s }-  b^{{}}c^{\s }] 
&\hspace{-6pt}+ &\hspace{-6pt}l [a^{\s}  d^{{} }- N   b^{\s }c^{{}}]\} .
\end{array}
\end{eqnarray}
 This presemifield multiplication can be generalized:
\begin{theorem}
\label{A''''theorem}
\label{Btheorem}
Let $L=\ef_{p^m}$ with $p$ odd,  $\s=p^s $ with $0< s<m, $ 
and     $l,n, N\in L^*, $ $ v\in L.$  Then
{\rm(\ref{genformequ})} and {\rm(\ref{A'''' mult})} 
define a presemifield $X(p,m,s,v,l,n,N)$
if and only if  
\begin{itemize} 
\item [\rm(i)]$-l \notin L^{\s-1} $ and
\item [\rm(ii)] the polynomial 
$ t ^{\s+{1}}-vt^{\s}-(v/N) t^{{}} +1/nN $
has no root   $t\in L$.
\end{itemize}
\end{theorem}

\begin{proof}
Assume (i), (ii), and that $(a,b)\star(c,d)=0$ with $(a,b),(c,d)\not=0.$ 
If $a=0,$ then $bc=-ad=0$ and $bd\not=0,$  in which case
$nN  b^{}d^{\s}+lnN b^{\s}d^{}=0$ is impossible, by hypothesis. 

If $a\not=0$ then we may assume that 
$a=1$ by homogeneity.  Since   $d=-bc,$ we obtain $c\ne0$ and
 $$
 \begin{array}{llll}
  \ \ \  \  \{[c^{\s }-n N  b^{{}}  (-bc)^{\s }] 
     &\hspace{-6pt}+\hspace{-6pt}&l[c^{{} }-  nN  b^{\s}  (-bc)^{{} }]\} \vspace{4pt}
     \\
 \hspace{-6pt} +nv \{[N  (-bc)^{\s }- \, b^{{}}c^{\s }] 
&\hspace{-6pt}+\hspace{-6pt}&l[ (-bc)^{{} }- N   b^{\s }c^{{}}]\}  =0.
\end{array}
$$
This simplifies to
$(c^{\s}+lc^{})\big (b^{\s+{1}}-vb^{\s}-(v/N) b^{{}} +1/nN \big) =0,$ which 
contradicts our   hypotheses.
The preceding argument reverses.
\end{proof}

\noindent{\bf Notation.}
Denote by $K_1=\ef_{p^{\gcd(m,s)}}=\ef_{q}\subseteq K_2=\ef_{p^{\gcd(m,2s)}}$ the fixed fields of
$\s$ and $\s^2$ in $L.$

\begin{remark}
Let $K\subset L$ be a field extension in odd characteristic.
We will repeatedly use the elementary fact that a non-square in $K$
is a non-square in $L$ provided $[L:K]$ is odd.
\end{remark}
 
\begin{Lemma}
\label{firstisolemma}
If $k_b,k_c\in L^* $ then the presemifields
$X(p,m,s,v,l,n,N)$ and  $X(p,m,s,v/k_b,l/k_c^{\s-1},nk_b^2,Nk_b^{\s-1})$ are isotopic$;$ 
and both are
isotopic to $X(p,m,m-s,v/N,1/l,nN^2,1/N).$
\end{Lemma}
\begin{proof} The first statement follows from the substitution
$$a'=a , \, b'=bk_b , \, c'=ck_c ,\,  d'=dk_d $$
$$l'=l/k_c^{\s -1}, \, v'=v/k_b, \, n'=nk_b^2,\,  N'=Nk_b^{\s -1}$$
where $k_d=k_bk_c,$ giving the product
$(h_{lvnN}(a',b',c',d'),a'd'+b'c')=(k_c^{\s}h_{l'v'n'N'}(a,b,c,d),k_bk_c(ad+bc)).$
(N.\:B.--The condition in Theorem~\ref{Btheorem}(ii) is unchanged: substitute $t \mapsto k_bt$ and divide
by $k_b^{\s +1}.$)
\par
Let $t=m-s, \tau =p^t$ and use the substitution 
$$a'=a^{\tau},\, b'=b^{\tau},\, c'=c^{\tau},\, d'=d^{\tau},$$
$$v'=v/N,\,  l'=1/l, \, n'=nN^2,\,  N'=1/N.$$
Then 
$$(h_{lvnN}(a',b',c',d'),a'd'+b'c')=(lh_{l'v'n'N'},(ad+bc)^{\s}).$$
This proved the second claim.
(This time the substitution $t \mapsto t^{\tau}$ shows that condition in Theorem~\ref{Btheorem}(ii) 
again is unchanged.)
\end{proof}

The special case $k_b=1$ shows in particular that $l$ can be replaced
by an arbitrary   element of  
$lL^*{}^{\s -1}$ without changing the isotopy type.
Lemma~\ref{firstisolemma} shows that all of the presemifields $X(p,m,s,v,l,n,N)$
in Theorem~\ref{A''''theorem}  are already isotopic to members 
of a significantly smaller subset:

\begin{Corollary} 
\label{BC corollary}
Every presemifield $X(p,m,s,v,l,n,N)$ is isotopic  to one of the following$:$

\begin{itemize} 
\item [\rm(i)] $B(p,m,s,l,n,N)=X(p,m,s,1,l,n,N),$ 
where  $0<s<m,$ $  -l\notin L^{\s -1}$ and
\begin{equation}
\label{like Knuth}
x^{\s +1}+(1-1/N)x+(1/n-1)/N \ne0 \ \mbox{ for  } \ x \in L;
\end{equation}
 or

\item [\rm(ii)] $C(p,m,s,l,R)=X(p,m,s,0,l,n,N),$ where
$0<s<m, -l\notin L^{\s -1}, R=-nN\notin L^{\s +1}$. 
\end{itemize}
\end{Corollary}
\begin{proof}
These follow from the lemma; in Theorem~\ref{A''''theorem}(ii) we have substituted 
  $t=x+1$  in order to obtain  the polynomial  in 
(\ref{like Knuth}). 
\end{proof}
 
\begin{Proposition}
\label{CBisoprop} 
  $C(p,m,s,l,R)$ is isotopic to $C(p,m,s,l',R')$  for any   $l'\in lL^*{}^{\s -1},$
  $R'\in RL^*{}^{\s +1}$.

$C(p,m,s,l,R)$ is isotopic to $C(p,m,m-s,1/l,R).$
\par
$B(p,m,s,l,n,N)$ is isotopic to $B(p,m,m-s,1/l,n,1/N^{1/\s}),$ and
  $l$ can be replaced by an arbitrary element of $lL^*{}^{\s -1}.$
\end{Proposition}
\begin{proof}
These follow from Lemma~\ref{firstisolemma}.
\end{proof}
 
\begin{remark} \rm
There is a curious similarity between the polynomials 
in $(\ref{like Knuth}) $
and  $(\ref{Knuth condition}),$ with $g=1-1/N $
and $f=(n-1)/nN $.    
However$,$ the present presemifields depend on more parameters in $L$ 
than Knuth's semifields.
\end{remark}

\section{The commutative case}
\label{commsection}

Recall that all commutative presemifields are isotopic to commutative ones
by (\ref{isotope identity})  or  (\ref{isotope identity e}).

\subsection{The commutative C-family}
\label{commCsub}

\begin{Theorem}
\label{Ccommtheorem}
$C(p,m,s,l,R)$ (where $-l\notin L^{\s-1},R\notin L^{\s+1}$) is isotopic to a commutative
semifield if and only if $l^{\s+1}R^{\s-1}\in L^{\s^2-1}.$
\end{Theorem}
\begin{proof}
Assume first that  $l^{\s+1}R^{\s-1}\in L^{\s^2-1}.$ Up to isotopy 
we may assume that $l^{\s+1}R^{\s-1}=1.$ (For, if 
$l^{\s+1}R^{\s-1}=z^{-(\s^2-1)} $ then we can use 
Proposition~\ref{CBisoprop} 
to replace $l$ by $lz^{\s-1}$, where 
$(lz^{\s-1})^{\s+1}R^{\s-1}=1$.)
Let $k=lR.$ Then $k^{\s}=R/l$, and the substitution 
$a\mapsto kb,$ $ b\mapsto a$
transforms the presemifield multiplication 
of $C(p,m,s,l,R)$ in  Corollary~\ref{BC corollary}(ii)  
and (\ref{genformequ}) into
$(h(kb,a,c,d),ac+kbd)$, where
$$h(kb,a,c,d)=Rad^{\s}+lRa^{\s}d+lRbc^{\s}+Rb^{\s}c.$$
Conversely, assume that $C(p,m,s,l,R)$ is isotopic to a commutative semifield.
Recall the presemifield multiplication for $C(p,m,s,l,R):$
$$(a,b)*(c,d)=(\underbrace{ac^{\s}+Rbd^{\s}+l\lbrack a^{\s}c+Rb^{\s}d\rbrack}_{h(a,b,c,d)} ,ad+bc).$$
We use Corollary~\ref{preGanleyprop}  with  $v=(v_0,v_1).$
As $x*1=(a+la^{\s},b),$ we have
$\a (x)=(a',b)$ where $a'+la'^{\s}=a.$ 
Further $v*x=(h(v_0,v_1,a,b),v_0b+v_1a)$ and
$\a (v*x)=(h'(v_0,v_1,a,b),v_0b+v_1a)$, where $h'$ is defined by 

\begin{equation}
\label{h'hequ}
h'(a,b,c,d)+lh'(a,b,c,d)^{\s}=h(a,b,c,d)
\end{equation}

The commutativity condition in Corollary~\ref{preGanleyprop} therefore becomes the two requirements
\begin{eqnarray}
\hspace{-27pt}h(h'(v_0,v_1,a,b),v_0b+v_1a,c,d)\hspace{-7pt}&=&\hspace{-7pt}h(h'(v_0,v_1,c,d),v_0d+v_1c,a,b)
\label{comm1} \\
\hspace{-27pt}h'(v_0,v_1,a,b)d+v_0bc+v_1ac\hspace{-7pt}&=&\hspace{-7pt}h'(v_0,v_1,c,d)b+v_0ad+v_1ac.
\label{comm2}
\end{eqnarray}

Consider (\ref{comm2}) as a function of $c,d$ alone. 
Then $h'(v_0,v_1,c,d)$ has the form
$h'(v_0,v_1,c,d)=f_{v_0,v_1}d+g_{v_0,v_1}c$, 
where we write $f_{v_0,v_1}=h'(v_0,v_1,a,b)/b-v_0a/b)$ and $g_{v_0,v_1}=v_0.$
Abbreviate this to $h'(c,d)=fd+g$ and use it in 
(\ref{h'hequ}):
$$
\begin{array}{lllll}h'(c,d)+lh'(c,d)^{\s}\hspace{-7pt}&=&\hspace{-7pt}lf^{\s}d^{\s}+fd+g+lg^{\s} 
\\
\hspace{-7pt}&=&\hspace{-7pt}h=Rv_1d^{\s}+lRv_1^{\s}d+v_0c^{\s}+lv_0^{\s}c. 
\end{array}$$
Comparing coefficients of $d$ and $d^{\s}:$ $Rv_1=lf^{\s}, $  $lRv_1^{\s}=f,$
which in case $v_1\not=0$ combine to $l^{\s+1}R^{\s-1}=1/v_1^{\s^2-1},$ as desired.
Assume $v_1=0.$ Then $v_0\not=0, f=0, h'(v_0,0,c,d)=v_0c.$ Then
 (\ref{comm2})  is  
satisfied, and  (\ref{comm1}) reads
$h(v_0a,v_0b,c,d)=h(v_0c,v_0d,a,b)$ for all $a,b,c,d,$ which is satisfied if and only if
$l=1/v_0^{\s-1}.$ In this case we may assume that $l=1.$ The condition $-1\notin L^{\s-1}$
implies $(p^m-1)/(p^{\gcd(m,s)}-1)$ odd, equivalently $m/\gcd(s,m)$ odd and therefore $L^{\s-1}=L^{\s^2-1}.$ It follows
$l^{\s+1}R^{\s-1}=R^{\s-1}\in L^{\s^2-1}$ as claimed.
\end{proof}

In particular, the $C$-family contains many semifields not isotopic to commutative ones.
A small example occurs with $p=3,m=4$ and $l=R$ of order $16.$
\medskip

\noindent{\bf Notation.}
Denote by $v_2(n)$ the highest power of $2$ dividing the integer $n.$  The following is elementary:

\begin{Lemma}
\label{msoddlemma}
For odd $p$ we have $\gcd(p^m+1,p^s+1)=2$ if $v_2(s)\not=v_2(m),$ whereas
$\gcd(p^m+1,p^s+1)=p^{\gcd(m,s)}+1$ if $v_2(s)=v_2(m).$ 
\end{Lemma}

\begin{Lemma}
\label{msodd2ndlemma}
Let $p$ be odd, $d=\gcd(s,m).$ 
\begin{itemize}
\item [\rm(i)]
If $v_2(d)=v_2(m),$ then $\gcd(2s,m)=d,~\gcd(p^m-1,p^s+1)=2$
and $L^{\s-1}\cap L^{\s+1}=L^{\s^2-1}.$
\item [\rm(i)]
If $v_2(d)<v_2(m),$ then $\gcd(2s,m)=2d,~\gcd(p^m-1,p^s+1)=p^d+1$
and $L^{\s^2-1}$ has index $2$ in $L^{\s-1}\cap L^{\s+1}.$
\end{itemize}
\end{Lemma}
\begin{proof} 
The statements concerning $\gcd(2s,m)$ are obvious. Let $g=\gcd(p^m-1,p^s+1).$ 
Then $g$ divides $\gcd(p^m-1,p^{2s}-1)=p^{\gcd(2s,m)}-1.$
If $v_2(d)=v_2(m),$ then $\gcd(2s,m)=d$ and $g\vert p^d-1\vert p^s-1.$ As $g$ also divides $p^s+1$ it follows $g=2.$
Let $v_2(d)<v_2(m).$ Then $v_2(d)=v_2(s)$ and $g\vert (p^{2d}-1),$ which divides $p^m-1.$
It follows $g=\gcd(p^{2d}-1,p^{s}+1).$ Lemma~\ref{msoddlemma} shows $p^d+1\vert g.$ Further
$\gcd(p^{d}-1,p^s+1)=2$ as $v_2(d)=v_2(s).$ This shows $g=p^d+1.$
\par
Clearly $L^{\s^2-1}\leq L^{\s-1}\cap L^{\s+1}.$
Let $a\in L^{\s-1}\cap L^{\s+1}.$ Then $a^{\s+1}\in L^{\s^2-1}$ and 
$a^{\s-1}\in L^{\s^2-1}.$
It follows $a^2\in L^{\s^2-1}.$ This shows $i=[L^{\s-1}\cap L^{\s+1}:L^{\s^2-1}]\leq 2.$  
Use the familiar formula $\vert L^{\s-1}\vert =(p^m-1)/\gcd(p^s-1,p^m-1)$ and analogous formulas
for the other cyclic goups involved. We see that $i=2$ if and only if
both of $v_2(\gcd(p^m-1,p^s-1))$ and $v_2(\gcd(p^m-1,p^s+1))$ are smaller than
$v_2(\gcd(p^m-1,p^{2s}-1)).$ By what we saw above this is satisfied if and only if $v_2(d)<v_2(m).$
\end{proof}

\begin{Lemma}
\label{sigmalemma}
Let $L=\ef_{p^m}$ with $p$ odd,  $n\in L$ a non-square$,$ $\s=p^s,$ $ 0\leq s<m$ and
$d=\gcd(s,m).$  
Then 
\begin{itemize}
\item [\rm(i)]
$-n^{(\s-1)/2}\in L^{\s-1}$ if and only if $s/d$ and $m/d$ odd;
and
\item [\rm(ii)]
$n^{(\s-1)/2}\in L^{\s-1}$ if and only if $s/d$ even ($m/d$ is then necessarily odd).
\end{itemize}
\end{Lemma}
\begin{proof} 
(i) The first statement means that there is an element $x\in L$ such that $(n/x^2)^{(\s-1)/2}=-1,$ or
equivalently that there is a non-square $n'\in L$ such that
$n'^{(\s-1)/2}=-1.$ This is equivalent to $v_2(p^s-1)=v_2(p^m-1)$ and to
both $s/d$ and $m/d$ being odd.
\par
(ii) The proof  is similar: an equivalent condition is
$v_2(\s -1)>v_2(p^m-1),$ which is equivalent to 
$$v_2(p^d-1)=v_2(p^m-1),~v_2(p^d-1)<v_2(\s-1);$$
and these are equivalent to $m/d$ odd and $s/d$ even, respectively.
\end{proof}

\begin{Proposition}
\label{Ccommprop}
Let $d=\gcd(s,m).$
The   $C$-family of commutative presemifields in
 {\rm Theorem~\ref{Ccommtheorem}}
splits into  two subfamilies$:$
\begin{itemize}
\item [\rm(i)]
$m/d$ is odd, $l=1$ and $R$ is a non-square in  
$K_1=\ef_{p^d},$ and
\item [\rm(ii)]
 $m/d$ is even and the presemifield is
$C(p,m,s,1/\o^{s-1},\o^{s+1}) $ with $s\leq m/2.$
\end{itemize}
\end{Proposition}
\begin{proof}
As mentioned in the proof of Theorem~\ref{Ccommtheorem}, it may be assumed
that $l^{\s+1}R^{\s-1}=1.$ 
Let $Z=\lbrace (l,R)\vert u(l,R):=l^{\s+1}=1/R^{\s-1},-l\notin L^{\s-1},R\notin L^{\s+1}\rbrace .$
Then $u(l,R)\in L^{\s-1}\cap L^{\s+1}.$ Recall from Proposition~\ref{CBisoprop} that
$l$ may be replaced by an arbitrary member of $lL^*{}^{\s -1},$ and $R$ may be replaced by an arbitrary
element of $RL^*{}^{\s +1}$ without changing the isotopy type.
Concretely, if $(l,R)\in Z,$ then $(lx^{\s-1},R/x^{\s+1})\in Z$ for
arbitrary $0\not=x\in L,$ $ u(lx^{\s-1},R/x^{\s+1})=u(l,R)x^{\s^2-1}$ and $C(p,m,s,lx^{\s-1},R/x^{\s+1})$ is isotopic to $C(p,m,s,l,R).$ It follows that $u(l,R)\in L^{\s-1}\cap L^{\s+1}$ can be chosen arbitrarily in its coset mod $L^{\s^2-1}.$ 
Now use Lemma~\ref{msodd2ndlemma}. 

If $v_2(d)=v_2(m)$ there is only one such coset and
we may assume $u=1.$ Then $l^{\s+1}=1,$ and $l$ is in the subgroup of order
$\gcd(p^m-1,\s+1)=2.$ As $l\not=-1$ it follows $l=1.$ 
The condition $-1\notin L^{\s-1}$ translates as $m/\gcd(s,m)$ odd, and this implies that
non-squares of $K_1$ are non-squares of $L.$ As $R^{\s-1}=1,$ it follows $R\in K_1.$
The second condition on $R$ is $R\notin L^{\s+1}.$ By Lemma~\ref{msodd2ndlemma} we have
$\vert L^{\s+1}\vert =(p^m-1)/\gcd(p^m-1,\s+1)=(p^m-1)/2,$ which means that $R$ is a non-square.
Both conditions on $R$ taken together are equivalent to $R$ being a non-square in $K_1.$
The choice of this non-square $R$ is irrelevant up to isotopy.

Now consider the case  $v_2(d)<v_2(m).$ There are two cosets of $L^{\s^2-1}$ in $L^{\s-1}\cap L^{\s+1},$
see Lemma~\ref{msodd2ndlemma}.
Assume  first that $u(l,R)=1;$
we will see that this cannot occur.  This time $l$ is in the group of order
$\gcd(p^m-1,\s+1)=p^d+1.$ We have $\gcd(p^m-1,p^{2s}-1)=p^{2d}-1.$
We can replace $l$ by $l'=lx^{\s-1}$ provided $x^{\s^2-1}=1.$
This means $x$ is in the group of order $p^{2d}-1$ and $x^{\s-1}$ in the group of order $p^d+1.$ We may therefore assume $l=1$, which forces $-1\notin L^{*\s-1},$ equivalently
$m/d$ odd, a contradiction to the current assumption.
\par
It follows that $u=u(l,R)$ is a representative of the nontrivial
coset of $L^*{}^{\s^2-1}$ in $L^*{}^{\s-1}\cap L^*{}^{\s+1}.$ We claim that
$u=1/\o^{\s^2-1}$ is a representative.
In fact, clearly $u\in L^*{}^{\s-1}\cap L^*{}^{\s+1}$ and $u^2\in L^{\s^2-1}.$ Assume $u\in L^{\s^2-1}.$ 
As $n=\o^2$ is a non-square in $L,$ Lemma~\ref{sigmalemma} implies $m/d$ odd which contradicts our
current assumption.
\par
We have $u=1/\o^{\s^2-1}.$ As $l^{\s+1}=u$ it follows $l=x/\o^{\s-1},$ where $x^{\s+1}=1.$
As $\gcd(p^m-1,\s+1)=p^d+1$ there are $p^d+1$ choices for $x.$ Up to isotopy we can replace
$l$ by $ly$ where $y\in L^{\s-1}$ and $y^{\s+1}=1.$ There are $p^d+1$ choices for $y.$ It follows
that we may assume $x=1,$ equivalently $l=1/\o^{\s-1}.$ The necessary condition $-l\notin L^{\s-1}$ is equivalent
to $m/d, s/d$ not both being odd. This is satisfied as $v_2(d)<v_2(m)$ (see Lemma~\ref{sigmalemma}).
Similarly $R^{\s-1}=1/u=\o^{\s^2-1}$ shows $R=\o^{\s+1}x$ where $x\in K_1.$
Observe that $K_1\subseteq L^{*\s+1}$ as follows from Lemma~\ref{msodd2ndlemma}. Those choices for $x$
therefore yield isotopic presemifields. We may choose $R=\o^{\s+1}.$
Assume the condition $\o^{\s+1}\notin L^{\s+1}$ is not satisfied: $(\o^2/x^2)^{(\s+1)/2}=1$ for some
$x\in L.$ This means there is a non-square $\nu =\o^2/x^2\in L$ such that $\nu^{(\s+1)/2}=1$ and therefore
$v_2(\s+1)>v_2(p^m-1).$ However $\gcd(p^m-1,\s+1)=p^d+1$ by Lemma~\ref{msodd2ndlemma}. It follows
$v_2(p^d+1)=v_2(p^m-1)$ which is not true as $p^{2d}-1$ divides $p^m-1.$
We have the uniquely determined semifields $C(p,m,s,1/\o^{\s-1},\o^{\s+1})$ where
$s\leq m/2, $ $v_2(d)<v_2(m).$ 
\end{proof} 
 
\noindent{\bf Notation.}
Denote the commutative semifields of Proposition~\ref{Ccommprop} by $C_{p,m,s}.$  There are
$\lfloor m/2\rfloor$ of them of   order $p^{2m}.$

\begin{Proposition}
\label{BuHeCprop}
The first family of commutative semifields from~{\rm\cite{BuHe11}}
coincides with the family of semifields $C_{p,m,s}$ of {\,\rm Proposition~\ref{Ccommprop}} satisfying
the additional  condition that $p\equiv 1$
{\rm mod~$4$}
if $m/d$ is odd and $s/m$ is even. 
\end{Proposition}

We note that the latter additional condition is not needed for the existence of these semifields.
\begin{proof}
In~\cite{BuHe11} it is shown that
$(bx)^{\s+1}-\ov{(bx)^{\s+1}}+\sum_{i=0}^{m-1}c_i(x\bar x)^{p^i}$
is a planar function on $F=\ef_{p^{2m}},$ where $b\in F^*, c_i\in L,$ $
y\to \sum_ic_iy^{p^i}$ describes an invertible linear mapping $L\to L$ and
$$\begin{array}{lll}
(i)\ \gcd(m+s,2m)=\gcd(m+s,m).
\vspace{2pt}
 \\(ii) \gcd(p^s+1,p^m+1)\not=\gcd(p^s+1,(p^m+1)/2).
\end{array}$$
The substitution $x\mapsto x/b$ shows that it can be assumed that $b=1.$
Observe that $x^{\s+1}-\ov{x^{\s+1}}\in L\o$, whereas
$\sum_ic_i(x\bar x)^{p^i} \in L.$
Applying the identity to $L\o$ and the inverse of the linear mapping to $L,$
we obtain the planar function $x^{\s+1}-\ov{x^{\s+1}}+N(x)$ where
$N\colon F\to L$ is the norm $N(x)=x\bar x$. After polarization
the commutative presemifield multiplication can therefore be chosen as
$x\circ y-\ov{x\circ y}+T(\ov{x}y),$
where $x\circ y=xy^\s + x^\s y$.  This is based on the compatibility
of $(\ov{x}y,L\o )$ and $(x\circ y,L).$ 
In $L$-coordinates the presemifield multiplication is
$$(Nad^{\s}+bc^{\s}+a^{\s}d+Nb^{\s}c ,ac-nbd)$$
where $n=\o^2, N=n^{(\s-1)/2}.$
\par
Condition (i) is equivalent to $v_2(s)\not=v_2(m),$
again equivalently: either $m/d$ is even or $s/d$ is even, and
$\gcd(\s+1,p^m+1)\not=\gcd(\s+1,(p^m+1)/2)$ is equivalent to
$v_2(\s+1)\geq v_2(p^m+1).$ 
\par
In terms of $L$-equations perform the substitution 
$a \mapsto b \mapsto -a/n$ and multiply the real part by $-n.$ 
The multiplication which results is
$$(ac^{\s}-nNbd^{\s}+(N/n^{\s-1})a^{\s}c-nb^{\s}d,ad+bc).$$
This is precisely the multiplication in $C(p,m,s,1/\o^{\s-1},-\o^{\s+1})=C_{p,m,s}.$ 
For any $s$ either itself or $m-s$ satisfies
conditions (i) and (ii). If $m/d$ is even, (ii) is automatically satisfied;
if  $m/d$ odd and $s/d$ even that condition is satisfied if and only if $p\equiv 1$ mod $4.$
\end{proof}

The smallest member of the commutative $C$-family that is not contained in the first family
of~\cite{BuHe11} has order $3^6.$ It is $C_{3,3,2}=C(3,3,2,1,-1)$ with presemifield 
multiplication
$$(a,b)*(c,d)=(ad^9+a^9d+bc^9+b^9c,ac-bd).$$

\begin{Proposition}
\label{Parthprop}
The family from~{\em \cite[Theorem~1]{JBprojsemi}} coincides with the
set of $C_{p,km',2k}$ with $m'>1$   odd.
\end{Proposition}
\begin{proof}
Let $\ef_q\subset L=\ef_{q^{m'}}\subset F=\ef_{q^{2m'}}$ 
for  odd $q,$ $m'>1.$ 
The family from~\cite[Theorem~1]{JBprojsemi}
relies on the compatibility of $(xy,L\o )$ and $(x\circ_Ay,L),$ 
using  $\s =q^2.$
The presemifield product therefore uses $\pi_\mathfrak{I}(x\circ_Ay)$ and $\pi_\mathfrak{R}(xy),$ 
hence can be written as
$$x*y=(Nad^{\s}+bc^{\s}+a^{\s}d+Nb^{\s}c,ac+nbd)$$
where $n=\o^2,N=\o^{\s-1}$ and $\o\in\ef_{q^2}.$
The substitution $a \mapsto b, b \mapsto a/n,$ followed by multiplying the
real part by $n$ transforms this into
$$(ac^{\s}+nNbd^{\s}+(1/N)(a^{\s}c+nNb^{\s}d),ad+bc)$$
which is the product of $C(p,km',2k,1/N,R)$ where $\s=q^2=p^{2k},$
$N=\o^{\s-1}$
and $R=nN=\o^{\s+1}.$ The condition $v_2(s)\not=v_2(m)$ is satisfied as
$s=2k$ and $m=km'$ and $m'$ odd. The second condition $\o^{\s+1}\notin L^{\s+1}$
is equivalent to $v_2(\s+1)\leq v_2(p^{km'}-1).$ This is satisfied as $v_2(\s+1)=2.$
We have $C(p,m'k,2k,1/N,R)$ where $m'>1$ is odd and $N=\o^{\s -1},R$ a non-square.
\end{proof}

It had been shown in~\cite{MP11} that the family from~\cite[Theorem~1]{JBprojsemi}
is isotopic to a subfamily of~\cite{BuHe11}.

\subsection{Commutative B-families}
\label{commBsub} 
 
 \begin{Theorem}
\label{Bcommtheorem}
$B(p,m,s,l,n,N)$ is isotopic to a commutative semifield if and only if
there is $v=(v_0,v_1)\not=0$ such that
the following equations hold using 
$f=nv_0+lnv_0^{\s}-lnNv_1^{\s}:$
\begin{eqnarray*}
lf^{\s}\hspace{-6pt}&=&\hspace{-6pt} nNv_0+lnNv_0^{\s}-nNv_1\\
v_0\hspace{-6pt}&=&\hspace{-6pt}lv_0^{\s}+nv_1-lnNv_1^{\s} .
\end{eqnarray*}
In particular, these hold when either
\begin{itemize}
\item [\rm(i)] $N^2=n^{\s-1},$ $ l\in (1/N)L^{\s-1},$ or
\item [\rm(ii)] $N=n^{\s-1},$ $ l\in L^{\s-1}.$
\end{itemize}
\end{Theorem}
\begin{proof}
We use Corollary~\ref{preGanleyprop}.
Here
$$x*1=(h(a,b,1,0),b)=(a+la^{\s}-n(b+lNb^{\s}),b).$$
Then $\a (x)=(a',b),$ where
$$a'+la'^{\s}=a+n(b+lNb^{\s}).$$
We have
$v*x=(h(v_0,v_1,a,b),v_0b+v_1a)$ and
$$\ \  \   \a (v*x)=(h'(a,b),v_0b+v_1a) , $$
where 
\begin{equation}
\label{h'hagainequ}
h'(a,b)+lh'(a,b)^{\s}=h(a,b)+n(v_0b+v_1a+lN(v_0b+v_1a)^{\s}).
\end{equation}

The left side of the  
condition in Corollary~\ref{preGanleyprop} is
$\a (v*x)*y=(h(h'(a,b),v_0b+v_1a,c,d),h'(a,b)d+(v_0b+v_1a)c).$
The right side is obtained by switching the roles of $x$ and $y.$
This leads to the following two conditions:
\begin{eqnarray}
\label{hh'firstequ}
h(h'(a,b),v_0b+v_1a,c,d)\hspace{-7pt}&=&\hspace{-7pt}h(h'(c,d),v_0d+v_1c,a,b)\\
\label{hh'secondequ}
h'(a,b)d+v_0bc\hspace{-7pt}&=&\hspace{-7pt}h'(c,d)b+v_0ad.
\end{eqnarray}

View the last equation as a function of $d$ alone. This shows
$$h'(c,d)=(h'(a,b)/b-v_0a/b)d+v_0c.$$
Here $f=f_{v_0,v_1}=h'(a,b)/b-v_0a/b$ is a constant, depending only on $v_0,v_1.$
Write $h'(c,d)=fd+v_0c, h'(a,b)=fb+v_0a$ and use~\ref{h'hagainequ}:
$$(fd+v_0c)+l(fd+v_0c)^{\s}=h(v_0,v_1,c,d)+n(v_0d+v_1c+lN(v_0d+v_1c)^{\s}).$$
Comparing coefficients of $d,d^{\s},c,c^{\s}$ produces $3$ equations, the coefficients of
$c,c^{\s}$ yielding the same equation. These are the relations given in the statement of the
theorem. 
Assume that they are satisfied.
Equation~\ref{hh'secondequ} is then satisfied. Equation~\ref{hh'secondequ} reads
$$h(fb+v_0a,v_0b+v_1a,c,d)=h(fd+v_0c,v_0d+v_1c,a,b)$$
Comparing coefficients this leads to $8$ equations, all of which are automatically satisfied.
\par
Let $v_0=0,v_1\not=0.$ Then $f=-nv_1, v_1^{\s-1}=1/(lN)$ (in particular $lN\in L^{\s-1}$).
Comparison with $lf^{\s}$ gives $N^2=n^{\s-1}.$

In case $v_1=0$ the condition is $N=n^{\s-1}$ and $l=1/v_0^{\s-1}.$
\end{proof}

\begin{Corollary}
\label{Bcommcor}
The exceptional cases of {\rm Theorem~\ref{Bcommtheorem}} lead to the following families of
commutative semifields:
\begin{itemize}
\item [\rm(1)] $B(p,m,s,1/\o^{\s-1},n,n^{(\s-1)/2}),$ where $n\in L$ a non-square
and either $m/\gcd(s,m)$ or $s/\gcd(s,m)$ even.
\item [\rm(2)] $B(p,m,s,-1/\o^{\s-1},n,-n^{(\s-1)/2}),$ where $n\in L$ a non-square,
$s/\gcd(s,m)$ odd.
\item [\rm(3)] $B(p,m,s,1,v^2,v^{\s-1})$ where $0\not=v\in L, m/\gcd(s,m)$ odd.
\item [\rm(4)] $B(p,m,s,1,n,n^{\s-1}),$ where $m/\gcd(s,m)$ odd.
\end{itemize}
\end{Corollary}
\begin{proof}
The first three cases subdivide Theorem~\ref{Bcommtheorem}(i).
Observe that the variant when $n$ is a square and $N=-n^{(\s-1)/2}$ does not occur
as in this case $-l\in L^{\s-1}.$
Recall that $l$ may be replaced by an arbitrary member of its coset $lL^*{}^{\s-1}$ (Proposition~\ref{CBisoprop}),
so that we may choose $l=1/N$ in the first three cases, $l=1$ in the last.
In the first case we may also choose $l=1/\o^{\s-1}$ or $l=\o^{\s-1},$ in the second
$l=-1/\o^{\s-1}$ or $l=-\o^{\s-1},$ and $l=1$ in the last two cases.
The conditions on $m/\gcd(s,m), s/\gcd(s,m)$ are equivalent to the condition
$-l\notin L^{\s-1}$ from Corollary~\ref{BC corollary}, see Lemma~\ref{sigmalemma}.
\end{proof}

Recall also that the members of the B-family need to satisfy the polynomial
condition in Corollary~\ref{BC corollary}.
Whenever $m/\gcd(s,m)$ is odd we can choose the non-square $\o^2\in K_1$ and then
$\o^{\s-1}=-1$ (see the remark preceding Lemma~\ref{firstisolemma}).
In the cases when $N^2=n^{\s-1}, l=1/N,$ a commutative presemifield is obtained by
the substitution $a \mapsto -nb,b \mapsto a$, which transforms the product to
$\big (h(-nb,a,c,d),ac-nbd \big)$, where
$$
\begin{array}{llll}
-h(-nb,a,c,d)=\\
\lbrack nNad^{\s}+na^{\s}d+nbc^{\s}+nNb^{\s}c\rbrack
+n\lbrace ac^{\s}+a^{\s}c+nNbd^{\s}+nNb^{\s}d\rbrace .
\end{array}
$$
\begin{Proposition}
\label{BuHeBprop}
The second family of commutative semifields from
{\rm\cite{BuHe11}} is contained in the union of the first family 
and the family from {\rm Corollary~\ref{Bcommcor}(1).}
\end{Proposition}

\begin{proof}
The PN function constructed in~\cite{BuHe11} is
$$\b x^{\s+1}+\ov{\b x^{\s+1}}+\g N(x)+\sum_{i=0}^{m-1}r_iN(x)^{p^i},
$$
where $\b\in F$ is a non-square, $\g\notin L,$ $r_i\in L$ and $\gcd(m+s,2m)=\gcd(m+s,m).$
Up to isotopy we can choose $\g =\o .$
Apply the inverse of the linear mapping $a+b\o\mapsto a+b\o +\sum_ir_ib^{p^i}.$
This produces the presemifield operation
$\b (x\circ y)+\ov{\b (x\circ y)}+\o T(\ov{x}y).$ This relies on the
compatibility of $(\ov{x}y,L\o )$ and $(\b (x\circ y),L\o ),$
or equivalently, $(\o xy,L)$ and  $(\o\b (\ov{x}\circ y),L).$ 
Here $x\circ y=xy^{\s}+x^{\s}y.$ The substitution $x\mapsto x/\o$ produces 
$A(p,m,s,1/\o^{\s-1},\b )$ (see Theorem~\ref{Atheorem}).
The condition $\gcd(m+s,2m)=\gcd(m+s,m)$ is equivalent to $v_2(s)\not=v_2(m),$
equivalently either $s/\gcd(s,m)$ or $m/\gcd(s,m)$ even and again equivalently
$-1/\o^{\s-1}\notin L^{\s-1}$ (see Lemma~\ref{sigmalemma}).
This verifies the first condition of Theorem~\ref{Atheorem}. The second condition is
satisfied as $\b$ is a non-square.
Because of isotopy we can write
$\b =vn+\o$ for some $v\in L,$ and our semifield is $X(p,m,s,v,1/\o^{\s-1},\o^2,\o^{\s-1}).$ 
The case $v=0$ leads to the first family as considered in Proposition~\ref{BuHeCprop}.
Let $v\not=0.$
 By Lemma~\ref{firstisolemma}, $X(p,m,s,v,1/\o^{\s-1},\o^2,\o^{\s-1})$
 is isotopic to
$$
\begin{array}{llll}
X(p,m,s,1,1/\o^{\s-1},\o^2v^2,\o^{\s-1}v^{\s-1})=\\
\hspace{1.5in}B(p,m,s,1/\o^{\s-1},\o^2v^2,\o^{\s-1}v^{\s-1}).
\end{array}
$$
\end{proof}

In particular, we obtain a parametric description of the second family from~\cite{BuHe11}.
In order to describe such a semifield of order $p^{2m}$ the following have to be chosen:

\begin{itemize}
\item A representative $s$ from $\lbrace s,-s\rbrace$ such that
either $m/d$ or $s/d$ is even, where $d=\gcd(s,m).$ Let $\s =p^s.$
\item An element $0\not=v\in L=\ef_{p^m}$ such that $\b=v\o^2+\o\in F$ is a non-square,
equivalently: $N(\b)=v^2\o^4-\o^2\in L$ is a non-square.
\item The semifield is $X(p,m,s,v,1/\o^{\s-1},\o^2,\o^{\s-1}),$ 
which is isotopic to
$B(p,m,s,1/\o^{\s-1},\o^2v^2,\o^{\s-1}v^{\s-1}).$
\end{itemize}

If $m/d$ is odd, we can choose $\o^2$ as a non-square in $\ef_{p^d}$
and obtain $\o^{\s-1}=-1.$
\par
We have not yet encountered an isotopy between 
semifields in
different cases in Corollary~\ref{Bcommcor}.
We  therefore expect commutative semifields of types 2, 3 and 4 
to be non-isotopic to members of the families from~\cite{BuHe11}.
\par
The semifield in Corollary~\ref{Bcommcor}(2) can be parametrized as\break 
$B(p,m,s,-1/\o^{\s-1},v^2\o^2,-v^{\s-1}\o^{s-1}).$ 
The polynomial condition (\ref{like Knuth}) reads  
$$t^{\s+1}+(v^{\s}+v^/\o^{\s-1})t+v^2/\o^{\s-1}-1/\o^{\s+1}\not=0\mbox{ for }t\in L.$$

\begin{Proposition}
Suppose that  $m/\gcd(s,m)$ and $ s/\gcd(s,m)$ are odd, $\o^2\in K_1$ is non-square and $v\in K_1$ is  such that
$v^2-1/\o^2$ is a non-square. Then the presemifield
$B(p,m,s,1,(v\o )^2,1)$ is defined and is isotopic to a commutative semifield.
\end{Proposition}
\begin{proof} This is contained in the family from Corollary~\ref{Bcommcor}(2).
We have $\o^2\in K_1$ and therefore $\o^{\s-1}=-1.$ The polynomial condition simplifies to
$v^2-1/\o^2\notin L^{\s+1}.$ As $\gcd(p^m-1,p^s+1)=2$ (see Lemma~\ref{msodd2ndlemma})
this is satisfied.
\end{proof}

Consider Corollary~\ref{Bcommcor}(3) where $n=v^2$ is a square. The polynomial condition 
(\ref{like Knuth})
reads
$$t^{\s+1}+(v^{\s}-v)t+1-v^2\not=0\mbox{ for }t\in L.$$
Computer experiments suggest that there are many values of $v$ which satisfy it.
An easy way to obtain examples is to choose $v\in K_1.$ The polynomial condition simplifies
then to $v^2-1\notin L^{\s+1}.$ 
It follows from Lemma~\ref{msodd2ndlemma} that
$\gcd(p^m-1,p^s+1)=2.$ The condition on $v$ is therefore: 
$v^2-1$ is a non-square in $K_1.$ 
\par
A similar procedure leads to the following parametric form of  Corollary~\ref{Bcommcor}(4):

\begin{Theorem}
\label{piuoltreBuHetheorem}
Let $m/d$ odd and $0\not=n\in L$ such that
$$t^{\s+1}+(n^{\s}-n)t+n-n^2\not=0\mbox{ for }t\in L.$$
Then $B(p,m,s,1,n,n^{\s-1})$ exists and is isotopic to a commutative semifield.
The above polynomial condition is satisfied in particular when $n\in K_1$ and $1-1/n\in K_1$ is a non-square.
\end{Theorem}

\section{The nuclei of $B(p,m,s,l,n,N)$   and\\  
$C(p,m,s,l,R)$}
\label{oddBCnucleisection}

First we handle a special case.
  
\subsection{The semifields $C(p,2k,k,l,R)$}

The semifields $C(p,2k,k,l,R)$ and $B(p,2k,k,l,n, N)$ of order $q^4$ (where $q=p^k$)
have $\s^2=1$   on $L.$
Here we  give a very  explicit description   of  the presemifields $C(p,2k,k,l,R)$, and use this to prove the  
\begin{proposition}
\label{special case of of nuclei}
 $\nz _l= \nz_r=\ef_q$ and $\nz_m=L=\ef_{q^2} $  for $C(p,2k,k,l,R).$
\end{proposition}

\proof
The conditions of Theorem~\ref{Btheorem} are: $N(l):=l^{q+1}\not=1$ and $R\notin  \ef_q.$  
Write $\bar a =a^\s =a^q$
for $a\in L$. 
By (\ref{B mult}), 
$$h(a,b,c,d)=\lbrack a \bar c +Rb \bar d \rbrack +l\lbrack \bar a c+R \bar b d\rbrack.  $$
 
Let $x=(a,b), y=(c,d), z=(e,f).$ Then
$1*x=(\bar a +la,b), x*1=(a+l \bar a ,b)$ and $\g (a,b)=(\bar a ,b).$
A direct calculation, using the fact that 
the inverse of $\bar x +lx$ on $L$ is 
$(\bar x -\bar l x)/(1-N(l) ),$ leads to the explicit 
{\em semifield} multiplication
$$x\circ y=(ac+R_1 \bar b d+R_2b \bar d , \bar a d+bc)  $$
with $R_1:=\frac{\bar R -N(l)R}{1-N(l)},$ $ R_2:=\frac{\bar l (\bar R -R)}{1-N(l)}.$  
Then 
$$\begin{array}{lllll}
x\circ (y\circ z)&\hspace{-6pt}=& \hspace{-6pt} \big(ace+R_1a \bar d f+R_2ad \bar f +R_1 \bar b (\bar c f+de)+R_2b(c \bar f + \bar d  \bar e ), 
\\ && \hfill \bar a (\bar c f+de)+b(ce+R_1 \bar d f+R_2d \bar f ) \big) \vspace{2pt}
\\
(x\circ y)\circ z&\hspace{-6pt}=& \hspace{-6pt} \big(ac+R_1 \bar b d+R_2b \bar d )e+R_1(a \bar d + \bar b  \bar c )f+R_2(\bar a d+bc) \bar f,  
\\&& \hfill(\bar a  \bar c + \bar R_1 b \bar d + \bar R_2  \bar b d)f+(\bar a d+bc)e \big).
\end{array}$$
Therefore $x\circ (y\circ z)=x\circ y)\circ z$  is equivalent to
$$(a-\bar a )d \bar f =b \bar d (e-\bar e ) ~\mbox{and} ~(R_1-\bar R_1 )b \bar d f= \bar R_2 \bar b df-R_2bd \bar f .$$
It follows easily that the nuclei are as stated.  \qed

\medskip
 Applying a suitable Knuth operation \cite[Sec.~4]{KnuthJAlg}
to $C(p,2k,k,l,R)$ produces a semifield with left nucleus of order $q^2$ that is 
a lifting of a Desarguesian plane  
 \cite[(8)]{CPT}.
 
\subsection{The generic case: preliminaries}

In odd characteristic,  define $(a,b)*(c,d)$ by {\rm(\ref{genformequ})} 
and {\rm(\ref{A'''' mult})}. 
This describes $B(p,m,s,l,n,N)$ (case $v=1$) and
$C(p,m,s,l,R)$ (case $v=0,$ where $R=-nN$).
We wish to determine the nuclei.
Recall that $\s\ne1$ on $L$.

Let $x=(a,b), y=(c,d), z=(e,f).$
We use the formulas
\begin{eqnarray*}
1*x=(h(1,0,a,b),b)=(a^{\s}+la+nv(Nb^{\s}+lb),b)\\
x*1=(h(a,b,1,0),b)=(a+la^{\s}-nv(b+lNb^{\s}),b) 
\end{eqnarray*}
in order to calculate the functions $\b,\g$ given in (\ref{beta gamma}).
Proposition~\ref{pretosemigenprop}  provides the semifield product
$x\circ y=\b (\g (x)*y)$,
where $\g (x)=(V,b) $ with 
\begin{equation}
\label{Vequ}
V+lV^{\s}-nv(b+lNb^{\s})=a^{\s}+la+nv(Nb^{\s}+lb) 
\end{equation}
and $\b (x)=(V',b)$ with
\begin{equation}
\label{V'equ}
V'^{\s}+lV'+nv(Nb^{\s}+lb)=a.
\end{equation}
Also $\g (y)=(W,d)$ with
\begin{equation}
\label{Wequ}
W+lW^{\s}-nv(d+lNd^{\s})=c^{\s}+lc+nv(Nd^{\s}+ld).
\end{equation} 
Lastly,  define $k$  and  $M$ by
\begin{eqnarray}
\label{kequ}
\hspace{-16pt}k+lk^{\s}-nv(Vd+bc+lN(Vd+bc)^{\s})\hspace{-6pt}&=&\hspace{-6pt}h(V,b,c,d)
\\
\label{Mequ}
\hspace{-16pt}M^{\s}+lM+nv(N(Wf+de)^{\s}+l(Wf+de))\hspace{-6pt}&=&\hspace{-6pt}h(W,d,e,f).
\end{eqnarray}
Note that all of the definitions (\ref{Vequ})--(\ref{Mequ})
use the fact that $x\to x+l x^{\s}$ is bijective on $L$ by 
Theorem~\ref{Btheorem}(i).
In terms of the notation  (\ref{Vequ})--(\ref{Mequ}) we have the
following
\begin{Lemma}
\label{assocbasiclemma}
The equation $x\circ (y\circ z)=(x\circ y)\circ z$
is equivalent to the following system$:$
\begin{eqnarray}
\label{Imequ}
bM \hspace{-6pt} &=& \hspace{-6pt}bce+(k-VW)f  \\
\label{Reequ}
h(k,Vd+bc,e,f)\hspace{-6pt} &=& \hspace{-6pt} h(V,b,M,Wf+de) .
\end{eqnarray}
\end{Lemma}
\begin{proof}
By (\ref{isotope identity e}), $x\circ (y\circ z)=(x\circ y)\circ z$ is equivalent to 
\begin{equation}
\label{seven4equ}
\g (x\circ y)*z=\g (x)*(y\circ z).
\end{equation}
 Since
$$\g (x)*y=(V,b)*(c,d)=(h(V,b,c,d) ,Vd+bc) $$
we have $x\circ y=(m,Vd+  b   c)$, where 
$m^{\s}+lm+nv(N(Vd+bc)^{\s}+l(Vd+bc))=h(V,b,c,d).$ 
Then
$\g (x\circ y)=(k,Vd+bc)$  by   (\ref{kequ}).

Also $\g (y)*z=(W,d)*(e,f)=(h(W,d,e,f) ,Wf+de),$ so that
 $y\circ z=(M,Wf+de)$ by   (\ref{Mequ}).
 
Now (\ref{seven4equ}) becomes
$(k,Vd+bc)*(e,f)=(V,b)*(M,Wf+de),$ and comparison of coordinates
completes the proof.%
\end{proof}

\begin{corollary}
\label{W and c}
If $x\circ (y\circ z)=(x\circ y)\circ z$ for some $x=(a,b),~b\ne0$, some $y$ and all 
$z$, then 
\begin{equation}
\label{Wnewequ}
\begin{array}{lll}
W\hspace{-6pt}&=&\hspace{-6pt} c^{\s}+nv(d+Nd^{\s}) \hspace{88pt}
\vspace{2pt}
\\
W^{\s}\hspace{-6pt}&=&\hspace{-6pt} c\hspace{5.7pt}+nv(d+Nd^{\s})
\end{array}
\end{equation}
\vspace{-8pt}
\begin{equation}
\begin{array}{lll}
\label{knewequ}
k\hspace{-6pt}&=&\hspace{-6pt} VW+nvb(c-c^{\s})-nNbd^{\s}
\vspace{2pt} 
\\ 
 k^{\s}\hspace{-6pt}&=&\hspace{-6pt} (VW)^{\s}+nvNb^{\s}(c^{\s}-c)-nNb^{\s}d .
\end{array}
\end{equation}
In particular$,$ \vspace{-6pt} 
\begin{eqnarray}
\label{c-c}
 c^{\s} -c\hspace{-6pt}&=&\hspace{-6pt}W-W^{\s} 
 \\
\label{cdfirstequ}
c-c^{\s^2}\hspace{-6pt}&=&\hspace{-6pt}(nv(d+Nd^{\s}))^{\s}-nv(d+Nd^{\s}) 
\end{eqnarray}
\vspace{-20pt}
\begin{equation}
\label{cdsecondequ}
(nv)^{\s}(c-c^{\s})^{\s}+nvN(c-c^{\s})=(nN)^{\s}d^{\s^2}-nNd.
\end{equation}
\end{corollary}
\begin{proof}  
As $b\not=0$,  (\ref{Imequ}) yields $M=ce+((k-VW)/b)f.$  
By definition, $V,W,k$ depend only on $a,b,c,d$ not $e$ or $f$. 
Then (\ref{Mequ}) yields 
(\ref{Wnewequ})--(\ref{knewequ})
 by comparing coefficients of $e,e^{\s},f,f^{\s}$ 
 and using the facts that $e,f\in L$ are arbitrary and $\s\ne1$
 on $L$.
 Now (\ref{c-c})--(\ref{cdsecondequ}) are immediate by (\ref{Wnewequ})--(\ref{knewequ}).
 \end{proof}

Later we will need the identities
\begin{eqnarray}
\label{emultequ}
\begin{array}{lll}
h(ea,eb,c,d)-eh(a,b,c,d)&  \vspace{3pt}
\\
\hspace{63pt}=l(e^{\s}-e)\lbrace (a^{\s}c-nNb^{\s}d)+nv(a^{\s}d-Nb^{\s}c)\rbrace 
\end{array} 
\\
\begin{array}{lll}
\label{emultequ2}
h(a,b,ec,ed)-eh(a,b,c,d)&  \vspace{3pt}
\\\hspace{63pt}=(e^{\s}-e)\lbrace (ac^{\s}-nNbd^{\s})+nv(Nad^{\s}-bc^{\s})\rbrace .
\end{array}
\end{eqnarray}
\vspace{-11pt}


\subsection{Left and right nuclei}

 By (\ref{B mult}),   
$(a,0)*(c,0)= (ac^\s+la^\s c,0)$, so there is a twisted field 
induced on $(L,0 )$.  {\em We identify the set 
$(L,0 )$ with $L$.}  By 
(\ref{beta gamma}) and 
(\ref{isotope identity e}), $(L,+,\circ )$ is a semifield associated with this twisted field.

Let $K_1=\ef_q$ and $K_2$ denote the respective fixed fields of $\s$ and $\s^2$ in $L$. 
 We frequently  identify these with $(K_1,0 )$ and $(K_2,0 )$, respectively.
These will be significant for calculating the nuclei.  First we note the 

\begin{Lemma}
\label{centerlemma}
$K_1$ is in the center of $(F,+,\circ ).$
\end{Lemma}
\begin{proof}
Apply Proposition~\ref{pretosemigenprop}(i).
\end{proof}

We will need the nuclei of the twisted field $L$
\cite[Theorem~1]{Alberttwist}: 
\begin{equation}
\begin{array}{lll}
\label{Albert}
\mbox{\em If $\s^2\ne1$ on $L$ then  $K_2 $ is the  middle nucleus of 
$(L,+,\circ )$}\\
\mbox{\em  while $K_1$ is both the left and right nucleus.}
\\
\mbox{\em  If $\s^2=1$ on $L$ then $(L,+,\circ )$ is a field.}
\end{array}
\end{equation}
\begin{Theorem}
\label{BClrnuketheorem}
The semifields    $B(p,m,s,l,n,N)$ 
and $C(p,m,s,l,R)$
have left and right nucleus $K_1\subset L.$ 
\end{Theorem}
\begin{proof} 
By Lemma~\ref{centerlemma}, $K_1$ is in both of these nuclei.\smallskip

{\noindent \bf Case $\s^2=1$ on $L.$}  
We need only consider $B(p,2k,k,l,n,N) $  in view of 
 Proposition~\ref{special case of of nuclei}. Recall $q=p^k.$
Let   $(a,b)\in \nz_l  $

First suppose that   $b\ne 0$.  Then   (\ref{cdfirstequ}) yields 
$n(d+Nd^q)=(n(d+Nd^q))^q$ for all $d,$ and hence $N=n^{q-1}.$  
Then $N=1$ by (\ref{cdsecondequ}), which contradicts Theorem~\ref{Btheorem}(ii). 
Thus   $\nz_l \subseteq L.$ 

Similarly,
$ \nz_r \subseteq L.$ 

\par
Assume that $(a,0)\in \nz_l.$ Then (\ref{Imequ}) gives $k=VW,$ 
while (\ref{Vequ}) gives 
 $V=  a ^q $  since $-l\notin L^{\s-1}$.  
Rewrite   (\ref{Mequ})
and   (\ref{Reequ}):    
\begin{eqnarray*}
V(M^q+lM)\hspace{-7pt}&=&\hspace{-7pt}Vh(W,d,e,f)-Vn\lbrace N(Wf+de)^q\!+l(Wf\!+de)\rbrace 
\\
VM^q\!+lV^qM\hspace{-7pt}&=&\hspace{-7pt}h(VW,Vd,e,f)\!-n\lbrace NV(Wf+de)^q\!+lV^q(Wf\!+de)\rbrace .
\end{eqnarray*}
{\em We claim that}  $V^q=V$ and hence $a^q=a$. 
 For if not,  subtract, use (\ref{emultequ}), factor out $l(V-V^q)$, and then use  (\ref{Mequ}) 
again  to get  
\begin{eqnarray*}
M\hspace{-6pt}&=&\hspace{-6pt} \lbrace W^q-n(d+Nd^q) \rbrace e-n\lbrace (W-W^q)+Nd^q\rbrace f \\
M^q\hspace{-6pt}&=&\hspace{-6pt}\lbrace W-n(d+Nd^q) \rbrace e^q-nN\lbrace (W^q-W)+d\rbrace f^q.
\end{eqnarray*}
Compare these:  the coefficient of $e^q$ shows $N=n^{q-1},$ 
and the coefficient of $f^q$ yields 
$n\in K_1=\ef_q.$ Then $N=1,$ which again contradicts Theorem~\ref{Btheorem}(ii).
\par
Finally, assume that $(e,0)\in \nz_r$.  {\em We claim that}   $e^q=e.$ 
For otherwise, (\ref{Imequ}) shows $M=ce.$ Comparison with (\ref{Mequ}) shows
$$We^q+lW^qe=(ce)^q+lce+n\lbrace de^q+N(de)^q+lde+lNd^qe\rbrace .
$$
Comparison with (\ref{Wequ}) 
produces
(\ref{Wnewequ}). 
This implies that $N=n^{q-1}.$
Rewrite (\ref{kequ}) and   (\ref{Reequ})  as 
\begin{eqnarray*}ek+elk^q \hspace{-6pt}&=&\hspace{-6pt}eh(V,b,c,d)+
en\lbrace (Vd+bc)+lN(Vd+bc)^q\rbrace 
\\
 ke^q+lk^qe\hspace{-6pt}&=&\hspace{-6pt}h(V,b,ce,de)+n\lbrace (Vd+bc)e^q+lN(Vd+bc)^qe\rbrace .
 \end{eqnarray*}
Subtract, use (\ref{emultequ2}) and divide by $e^q-e$: $$k=Vc^q-nNbd^q+n(NVd^q+Vd-bc^q+bc).$$
Then (\ref{kequ})  implies that
$$k^q=V^qc-nNb^qd+n(V^qd-Nb^qc+N(Vd+bc)^q\rbrace .$$
Comparison yields the usual contradiction $N=n^{q-1}=1,$
which finally proves our claim.
\smallskip

{\noindent \bf Case $\s^2\ne 1$ on $L.$} 
By Lemma~\ref{centerlemma} and (\ref{Albert}), 
$\nz_l\cap L =\nz_r\cap L = K_1.$ 
We must  show that $\nz_l$ and $\nz_r$ are in $L.$

Assume that $(a,b)\in \nz_l$ with $b\not=0.$  
Choosing $d=0$ in (\ref{cdfirstequ}) shows that $c\in K_2$.
Since  $c$ is arbitrary, we have $K_2=L $ and hence $\s^2=1 $
 on $L,$ which is not the case.  Thus, $\nz_l\subseteq (L,0)$, as required.
 
Assume that $(e,f)\in \nz_r$ with $ f\not=0.$
Since $\g(x)=(V,b)$, $c,d,b,V\in L$ are arbitrary.  Then  (\ref{Imequ}) yields
$k=VW+\big((M-ce)/f \big)b.$ 
As in the proof of Lemma~\ref{W and c},  
comparison of the coefficients of  $V,V^{\s},b$ and $b^{\s}$ in 
(\ref{kequ}) yields four equations.  Two of these are
(\ref{Wnewequ}), which we compare and choose $d=0$ in order to
deduce that $\s^2=1$ on $L$, which is not the case.
\end{proof}

\subsection{Middle nucleus}
Our results for the middle nucleus are not as satisfactory as for the left and right nuclei.
We start by identifying the middle nucleus in terms of field equations.

\begin{Lemma}
\label{cdNmlemma}
$(c,d)\in \nz_m$ if and only if 
{\rm(\ref{cdfirstequ})} and {\rm(\ref{cdsecondequ})}
hold.
\end{Lemma}
\begin{proof}
Let $(c,d)\in \nz_m .$ Choosing $b\not=0 $ yields 
 (\ref{Wnewequ})--(\ref{cdsecondequ}).
\par
Conversely, assume that {\rm(\ref{cdfirstequ})} and {\rm(\ref{cdsecondequ})} hold.
If $\tilde W:=c^{\s}+nv(d+Nd^{\s})$, then (\ref{cdfirstequ}) implies that both parts of
(\ref{Wnewequ}) hold with $W$ replaced by $\tilde W$.  Then 
(\ref{Wequ}) holds with $W$ replaced by $\tilde W$, so the uniqueness of $W$ 
in the definition (\ref{Wequ}) implies that
$\tilde W =W$: (\ref{Wnewequ}) holds.  The same reasoning yields 
(\ref{knewequ}).

In order to show that $(c,d)\in\nz_m$ we need to prove that 
{\rm(\ref{Imequ})} and {\rm(\ref{Reequ})} are satisfied for all $V,b,e,f.$ 

Assume at first $b\not=0.$ 
Let $\tilde M:= ce+ \big((k-VW)/b \big)f=ce+ \big(nv(c-c^{\s})-nNd^{\s} \big)f $    
(using  (\ref{knewequ})).  Then (\ref{Mequ}) holds with $\tilde M$ in place of $M$, and hence 
$\tilde M = M$ and (\ref{Imequ})  is satisfied.
We now have   $W,W^{\s},k,k^{\s},M$, and a calculation produces
{\rm(\ref{Reequ})}. (In fact, $h(k,Vd+bc,e,f)-h(V,b,M,Wf+de)$ is a linear
combination of $e,e^{\s},f,f^{\s}$ with coefficients depending only on $V,b,c,d,$ and 
a straightforward  calculation shows that all four coefficients vanish.)

Finally, let $b=0.$ Then $k=VW$  by (\ref{Wequ}), so that (\ref{Imequ}) holds.
As above, 
$M=ce+\big(nv(c-c^{\s})-nNd^{\s} \big)f$ follows from {\rm(\ref{Mequ})}, and then 
 (\ref{Reequ})  follows easily.
\end{proof}
 
\begin{Theorem}
\label{BspecialNmtheorem}
\begin{itemize}
\item [\rm(i)]
For $C(p,m,s,l,R),$ $\nz_m$ is $ K_2$ if $R\notin K_1^*L^{*\s +1},$ 
and a quadratic extension of $K _2$ otherwise.
\item [\rm(ii)]
For
$B(p,m,s,l,n,n^{\s -1}),$   $\nz_m$ is $ K_2$ if $ 1-1/n \notin K_1^{*}L^{*\s +1} $
and a quadratic extension of $K _2$ otherwise.
\item [\rm(iii)] For $B(p,m,s,l,n,N)$ with $N\not=n^{\s -1},$  
$\nz_m\cap L=K_1.$ 
\item [\rm(iv)]
For
$B(p,2k,k,l,n,N) $  with  
$N\ne n^{\s -1} = n^{q -1} ,$  $\nz_m$ is  a quadratic extension of $K _1 =\ef_q$.

\end{itemize}
\end{Theorem}
\begin{proof}
(i)  Here 
$v=0$, $R=-nN\notin L^{\s +1} $  by
Corollary~\ref{BC corollary}(ii).
Then (\ref{cdfirstequ}) yields $c\in K_2,$ 
hence $\nz_m\cap L \subseteq K_2.$

Assume that $d\not=0.$
Then (\ref{cdsecondequ}) yields
$d^{\s^2-1}=1/(nN)^{\s -1}=1/R^{\s -1},$ or equivalently $d^{\s +1}R\in K_1^*.$ 
There is a solution
$d\in L^*$ if and only if $R\in K_1^*L^{*\s +1}.$ 
Any two solutions $d_1,d_2$ satisfy $(d_1/d_2)^{\s^2-1}=1$ and hence 
$d_1/d_2\in K_2$. Thus, $|\nz_m|=|K_2|^2$ and $\nz_m$ has basis
$\lbrace (1,0),(0,d_0)\rbrace$ over $K_2,$ where $d_0$ is such that
$R\in d_0^{\s+1}K_1.$

(ii) Let $v=1,nN=n^{\s}.$ 
If $d=0 $, then  (\ref{cdfirstequ}) yields  $c\in K_2$ and 
(\ref{cdsecondequ}) implies $\nz_m\cap L=K_2.$ 

Suppose that  $d\not=0.$ Then (\ref{cdfirstequ}) simplifies to
$c-c^{\s^2}=(nd)^{\s^2}-nd,$ equivalently $c=-nd+t$ where $t\in K_2 \subseteq \nz_m .$
We may thus assume that  $c=-nd.$ 
Now (\ref{cdsecondequ}) states that $n^{\s}(c-c^{\s})=(nd)^{\s^2}-n^{\s}d.$ Using $c=-nd\ne 0$,
this simplifies to
$(dn)^{\s^2-1}=(1-1/n)^{1-\s},$ or equivalently $(dn)^{\s +1}(1-1/n)\in K_1.$
A solution $d\in L^*$ exists if and only if
$ 1-1/n \in K_1^{*}L^{*\s +1},$ and any two solutions $d_1,d_2$ satisfy $d_1/d_2\in K_2.$
Once again this implies that $\nz_m$ is a quadratic extension of $K_2.$

(iii)
Use Lemma~\ref{cdNmlemma} with $d=0.$

(iv)
 Let $(c,d)\in \nz_m  $.   Write $\bar a =a^\s =a^q$
for $a\in L$. 
Conditions 
(\ref{cdfirstequ})  and  (\ref{cdsecondequ}) state that 
$ n(d+N\bar d )\in \ef_q$ and 
$\bar n(\bar c-c)+nN(c-\bar c)= \bar n \bar N  d-nNd$. 
If $d=0$ then $c\in \ef_q$.  

Suppose that $d\ne 0$.  Then
$\bar d/d=(\bar n-nN )^{q-1} $, hence $d=(\bar n-nN )t, t\in K_1.$
Also, $\bar c-c = (\bar n\bar N-nN)t$ and $c=nNs, s\in K_1.$
It follows that $\nz_m$ is $2$-dimensional over $K_1,$ with basis  
$\lbrace (1,0),(nN,\bar n-nN)\rbrace .$
\end{proof}


We now assume that $\s^2\ne  1$ on $L$ and, until Proposition~\ref{last case}, 
also that  $  n^{\s-1}\ne   N^2$.  
Let  $\wz$ denote the kernel of the 
\vspace {3pt}
$K_1$-linear map $t\to (1+n\a_0)t+(1+n\a_1+nN\a_0^{\s})t^{\s}+(n\a_2+nN\a_1^{\s})t^{\s^2}+
nN\a_2^{\s}t^{\s^3},$  where
\begin{equation}
\label{alphacoeffequ}
\begin{array}{lllll}
\displaystyle\vspace{4pt}
\a_0\hspace{-6pt}&=&\hspace{-6pt}
\frac{\displaystyle n^{\s}-n^2N^2}
{\displaystyle n(nN^2-n^{\s})}\\
\vspace{4pt}
\displaystyle
~\a_1\hspace{-6pt}&=&\hspace{-6pt}
\frac{\displaystyle n^{\s}+n^2N^2-nN-n^{\s+1}N}
{\displaystyle n(nN^2-n^{\s})}\\
\displaystyle
\a_2\hspace{-6pt}&=&\hspace{-6pt}
\frac{\displaystyle nN(n^{\s}-1)}
{\displaystyle n(nN^2-n^{\s})} \,.
\end{array}
\end{equation}

\begin{Theorem}
\label {not N squared}
$\dim _{K_1}\nz_m={1+\dim _{K_1}\wz\leq 4}$
if 
$  n^{\s-1}\ne N, N^2$
and $\s^2\ne 1.$
\end{Theorem}

\begin{proof} 

First note that (\ref{cdfirstequ})
is equivalent to
$c+c^{\s}+n(d+Nd^{\s})\in K_1.$
For, if  $\a :=n(d+Nd^{\s}) $ then $\a^{\s}-\a =c-c^{\s^2} $ 
and hence $\a^{\s^2}-\a =(c-c^{\s})-(c-c^{\s})^{\s^2}$.
Then  $e:=\a +c+c^{\s} \in K_2.$ Substituting $\a =-c-c^{\s}+e$ in (\ref{cdfirstequ})
yields $e\in K_1.$

Since   $\nz_m\cap L=(K_1,0) $ by Lemma~\ref{centerlemma} 
and (\ref{Albert}) (as $\s^2\ne1$ on $L$), we may assume that $e=0$. Thus,  
 we are left to consider the $K_1$-space $\mz$ of pairs 
 $(c,d)\in L^2$ satisfying
\begin{eqnarray} 
\label{dsigmsequ1}
d^{\s}\hspace{-6pt}&=&\hspace{-6pt}-d/N-(c+c^{\s})/(nN)
\\\label{dsigmsequ2}
(nN)^{\s}d^{\s^2}\hspace{-6pt}&=&\hspace{-6pt}nNd+n^{\s}(c-c^{\s})^{\s}+nN(c-c^{\s}).
\end{eqnarray}
It is easy to check that  $\mz\cap (K_1,0)=0$, so that  $\dim \nz_m=\dim \mz+1$.  
Therefore, {\em it remains to prove that  
$\dim _{K_1} \mz=\dim _{K_1} \wz$.}

If (\ref{dsigmsequ1}) holds then  (\ref{dsigmsequ2})
is equivalent to  
\begin{eqnarray}
\label{d}
\begin{array}{llll}
n (nN^2-n^{\s})d =\big(n^{\s}- n^2N ^2 \big)c \\
\hspace{35pt}
+(n^{\s}+ n^2N ^2-nN-n^{\s+1}N)c^{\s}
+nN(n^{\s}-1)c^{\s^2}. 
\end{array}
\end{eqnarray}
Then $d=\a_0c+\a_1c^{\s}+\a_2c^{\s^2} $ with $\a_i$ in (\ref{alphacoeffequ}).
Substituting this into (\ref{dsigmsequ1}) shows that 
  $c\in \mz$ whenever $(c,d)\in \wz$, and this argument reverses. 
\end{proof}
\medskip

Note that we could have allowed $n^{s-1}= N$ in the above argument, but that case was already 
handled in Theorem~\ref{BspecialNmtheorem}(ii).
 
\begin{proposition} 
\label{last case}
If 
$  n^{\s-1}\!= N^2$
and $\s^2\ne 1 ,$ then $1\le \dim _{K_1}\!\nz_m\le3.$
\end{proposition}
\begin{proof} 
This time the left side of (\ref{d}) is 0, so that   (\ref{d}) describes a space  of solutions 
$c$ of $K_1$-dimension at most $2.$  
For each $c$ there are 0 or $|K_1|$ possible $d$ in  (\ref{dsigmsequ1}).
Hence, even allowing $d$ to be 0, the number of possible $(c,d)$ is at most $|K_1|^3$.
\end{proof} 

Recall from
the Motivation at the start of 
 Section~\ref{Generalization for all odd $q$} that 
$N=n^{(\s-1)/2} $   was one of the assumptions 
in Theorem~\ref{Atheorem}.

\end{document}